\documentclass[11pt,reqno,oneside]{amsart}
\usepackage{amscd}
\usepackage{amsmath}
\usepackage{latexsym}
\usepackage{amsfonts}
\usepackage{amssymb}
\usepackage{amsthm}
\usepackage{graphicx}
\usepackage{verbatim}
\usepackage{mathrsfs}
\usepackage{enumerate}
\usepackage{hyperref}
\usepackage{fullpage}
\usepackage{xcolor}
\usepackage[latin1]{inputenc}

\theoremstyle{plain}
\theoremstyle{plain}\newtheorem{theorem}{Theorem}[section]
\theoremstyle{plain}\newtheorem{lemma}[theorem]{Lemma}
\theoremstyle{plain}\newtheorem{coro}[theorem]{Corollary}
\theoremstyle{plain}\newtheorem{proposition}[theorem]{Proposition}
\theoremstyle{plain}\newtheorem{remark}{Remark}[section]

\numberwithin{equation}{section}

\newcommand{\R}{\mathbb{R}}
\newcommand{\be}{\begin{equation}}
\newcommand{\ee}{\end{equation}}
 \newcommand{\ba}{\begin{aligned}}
 \newcommand{\ea}{\end{aligned}}

  \newcommand{\ben}{\begin{enumerate}}
   \newcommand{\een}{\end{enumerate}}

\newcommand{\Rmnum}[1]{\expandafter\@slowromancap\romannumeral #1@}
\allowdisplaybreaks

\begin{document}
\title{Global well-posedness to the two-dimensional incompressible vorticity equation in the half plane}
\author{Quansen Jiu$^{1}$,\,You Li$^{2}$,\,Wanwan Zhang$^{3}$}

\address{$^1$ School of Mathematical Sciences, Capital Normal University, Beijing, 100048, P.R.China}
\email{jiuqs@cnu.edu.cn}

\address{$^2$ Department of Mathematics, Beijing Technology and Business University, Beijing, 100048, P.R.China}
\email{thu3141@126.com}

\address{$^3$School of Mathematical Sciences, Capital Normal University, Beijing, 100048, P.R.China}
\email{zhangwanwan153@163.com}

\begin{abstract}
 This paper is concerned with  the global well-posedness of the two-dimensional  incompressible vorticity equation in the half plane. Under the assumption that the initial vorticity $\omega_0\in  W^{k,p}(\R^{2}_+)$ with $k\geq3$ and $1<p<2$, it is shown that the two-dimensional incompressible vorticity  equation admits a unique solution $\omega\in C([0,T];W^{k,p}(\R^{2}_+))$ for any $T>0$. An elementary and self-contained proof is presented and delicate estimates of the velocity and its derivatives are obtained in this paper. It should be emphasized that the uniform estimate on $\int^t_0\|u(\tau)\|_{W^{1,\infty}(\R^2_+)}d\tau$ is  required to complete the global regularity of the solution.  To do that,  the double exponential growth in time of the gradient of the vorticity in the half plane is established and applied. This is different from  the proof of global well-posedness of the Euler velocity equations in the Sobolev spaces, in which a  Kato-type or logarithmic-type estimate of the gradient of the velocity is enough to close the energy estimates.
\end{abstract}
\smallskip
\maketitle

\section{Introduction and Main Results}
The  two-dimensional incompressible Euler equations read as
\begin{equation}\label{1.1}
\left\{\ba
&\partial_t u+(u\cdot\nabla)u+\nabla p=0, ~(x,t)\in D\times \R_+,\\&{\rm div} ~u=0,\ea\ \right.
\end{equation}
where  $D\subseteq \R^2$  is a domain and the unknown functions are the velocity field $u(x,t)= (u_1(x,t), u_2(x,t))$ and the pressure function $p$. The second equation means that the flow is incompressible, which enables us to determine the pressure $p$ from $u$ through a singular integral operator for a specific domain $D$. The initial condition to \eqref{1.1} is imposed as
 \begin{eqnarray}\label{1.2}
 u(x,0)=u_0(x),~x\in D.
 \end{eqnarray}
 Furthermore, when $D$ is a domain with a boundary, the natural boundary condition is the no  penetration one:
\begin{eqnarray}\label{1.3}
u(x,t)\cdot \nu(x)=0,~(x,t)\in\partial D\times \bar{\R}_+.
\end{eqnarray}
Here $\nu(x)$ denotes the outward unit normal vector of  domain $D$ at $x\in\partial D$.

 In two-dimensional case, the vorticity $\omega=\partial_{x_2}u_1-\partial_{x_1}u_2$ is a scalar function satisfying
\begin{equation}\label{1.4}
\left\{\ba
&\partial_t\omega+u\cdot\nabla\omega = 0, ~(x,t)\in D\times\R_{+}\\&u =\nabla^\bot\Delta^{-1}_D\omega,\\
&\omega(x,0)=\omega_{0}(x). \ea\ \right.
\end{equation}
The second equation of \eqref{1.4} is the so-called Biot-Savart law with $\nabla^\bot=(\partial_2,-\partial_1).$
For smooth domains $D$ with boundaries, $\Delta^{-1}_D$ denotes the inverse of the  Laplacian operator with the homogeneous boundary condition.
It is deduced from \eqref{1.4} that the quantity  $\|\omega(t)\|_{L^s(D)} (1\le s\le \infty)$ is conserved for all times, which plays a vital role in the proof of global regularity of solutions to the two-dimensional incompressible Euler equations.

 The global well-posedness theory to the Cauchy problem \eqref{1.1}-\eqref{1.2} or initial-boundary problem \eqref{1.1}-\eqref{1.3} has been established  in various settings.
 H\"{o}lder \cite{[Ho]} and Wolibner \cite{[W]} independently obtained the global existence and uniqueness of classical solutions to \eqref{1.1}-\eqref{1.3} in H\"{o}lder spaces for a smooth bounded domain $D$.
The global well-posedness result in \cite{[W]}  was later improved by Kato in \cite{[Kato1]} without the zero-circulation assumption for each inner component of the boundary of a bounded (not necessarily simply connected) domain, where the author used the Schauder fixed-point theorem and constructed the solutions   through the vorticity equation.
In \cite{[Mc]}, McGrath utilized the Schauder fixed-point theorem as well, but considered the vorticity equation in the stream function form to prove the  global well-posedness of classical solutions to \eqref{1.1}-\eqref{1.2} with $D=\R^2$.
Later, Bourguignon and Brezis \cite{[BB]} proved the global existence and uniqueness of strong solutions to \eqref{1.1}-\eqref{1.3} with $D$ a smooth bounded domain in Sobolev spaces with a suitable external force for the  initial velocity in $W^{m,s}(D)$ with $1<s<\infty$ and $m>1+\frac{s}{2}$.
For a recent proof of the global regularity, we also refer the reader to \cite{[Koch]}, where Koch proved the global existence and uniqueness of $C^1$ classical solutions to  \eqref{1.1}-\eqref{1.3} with a smooth bounded domain $D$ for $C^{1,\alpha}$ initial velocity  via the method of Schauder fixed-point theorem.
Existence and uniqueness of smooth solutions to \eqref{1.1}-\eqref{1.3} on exterior domain $D$ were obtained by Kikuchi in \cite{[Kik]}.
Concerning the half plane $D=\R^2_+$,  Secchi \cite{[S]} showed the existence and uniqueness of strong solutions to \eqref{1.1}-\eqref{1.3}  via the method of vanishing viscosity limit under specific Navier-boundary conditions.
In addition, it is well-known that Yudovich \cite{[Y]} obtained the existence and uniqueness of global weak solutions if the initial vorticity $\omega_0$ lies in $L^1\cap L^\infty$ for the domain $D=\R^2$ (see \cite{[MP]} for the bounded domain case). One can refer to
\cite{[Bar],[BC],[BL],[Che],[DM],[EJ],[EM],[HR],[KS],[KZ],[LZ],[X],[Z]} and the references therein for other related interesting and important aspects concerning with the two-dimensional incompressible Euler equations.

  Our goal in this paper is to consider the global well-posedness of strong solutions to the two-dimensional Euler in the half plane via the vorticity equation \eqref{1.4}. As we know, in the half plane, it has been  widely accepted and used in the boundary layer theory that there exists a unique global strong (classical) solution to the two-dimensional incompressible Euler equations (see e.g. \cite{[M1]} and \cite{[M2]}). And as mentioned above, it can be proved by using the  viscous approximation (the Navier-Stokes equations) with Navier boundary condition (see \cite{[L]} and \cite{[S]}) to the Euler equations. However, we will utilize the vorticity equation instead of the velocity equation to present an elementary and self-contained proof the  global well-posedness in this paper. Note that the vorticity equation is  one of a class of transport equations:
\begin{equation}\label{1.5}
\left\{\ba
&\partial_t\omega+u\cdot\nabla \omega = 0, ~(x,t)\in D\times\R_{+},\\&u =\nabla^\bot\Delta_D^{-1+\alpha}\omega,\\
&\omega(x,0)=\omega_{0}(x), \ea\ \right.
\end{equation}
where $0\leq\alpha\le\frac12$ and the domain $D$ is the whole plane, the torus, the upper half plane or the bounded domain. In fact, when $\alpha=0$, \eqref{1.5} is the two-dimensional incompressible vorticity equation. When  $\alpha=\frac12$, \eqref{1.5}   reduces to the celebrated two-dimensional inviscid  $\text{SQG}$  equation.
When $0<\alpha<\frac{1}{2}$, \eqref{1.5} is the so-called  inviscid modified or generalized $\text{SQG}$ in the literature (see \cite{[KYZ]} and references therein). There have been a number of  mathematical studies on inviscid SQG and inviscid modified SQG equations and we refer the readers to
\cite{[CC],[CCG],[CIN],[CMT],[CN1],[CN2],[C],[CF],[EJ],[GP],[HK],[I],[JK],[KN],[KRYZ],[KYZ],[M],[R]}
and the references therein for more details.

Our main results are stated as follows.

\begin{theorem}\label{the-1}
For every $\omega_{0}\in W^{k,p}(\R^{2}_+)$ with $k\geq3$ and $1<p<2$, there exists a  time $T_0=T_0(\|\omega_{0}\|_{W^{k,p}(\R^{2}_+)})$ such that \eqref{1.4} admits a unique solution $\omega\in C([0,T_0];W^{k,p}(\R^{2}_+))$.
\end{theorem}
\begin{theorem}\label{the-2}
For every $\omega_{0}\in W^{k,p}(\R^{2}_+)$ with $k\geq3$ and $1<p<2$. For any $T>0$, \eqref{1.4} admits a unique solution $\omega\in C([0,T];W^{k,p}(\R^{2}_+))$.
\end{theorem}
\begin{remark}
The restriction of $p\in(1,2)$ in Theorems \ref{the-1} and \ref{the-2} is required when estimating  $\|u\|_{L^{\frac{2p}{2-p}}(\R^{2}_+)}$ in terms of $\|\omega\|_{L^{\frac{2p}{2-p}}(\R^{2}_+)}$ in Proposition \ref{Sobolev-type estimate} with help of the Hardy-Littlewood-Sobolev inequality and when establishing the maximum estimate of the velocity and its gradient in Lemma \ref{Kato inequality}.
It would be interesting to obtain the local (global) well-posedness to \eqref{1.4} in the Sobolev spaces $W^{k,p}(\R^{2}_+)$ for more general values of $k$ and $p$, especially for $k\geq3$ and $2\leq p<\infty$.
\end{remark}
\begin{remark}
Theorems \ref{the-1} and \ref{the-2} can be also established for \eqref{1.4} with $D=\R^2$ by the same but simpler arguments .
\end{remark}
Now we explain our approaches in a more detail. To prove Theorem \ref{the-1}, starting from  the linearized equation \eqref{4.2},  we construct a suitable  approximate system (see \eqref{A.3})  and obtain the uniform estimate of the approximate solutions. In particular, we will establish a delicate Sobolev-type estimate on the singular integral involving the expression $u =\nabla^\bot\Delta^{-1}_D\omega$ in \eqref{1.4} with $D=\R^2_+$, that is,
\begin{eqnarray*}
\|u\|_{W^{m,\frac{2p}{2-p}}(\R^2_+)}\leq C_{m,p}\|\omega\|_{W^{m,p}(\R^2_+)},
\end{eqnarray*}
for any integers $m\geq0$ and $1<p<2$ (see Proposition \ref{Sobolev-type estimate}).
With the global well-posedness result of the linearized equation in hand (see Lemma \ref{Linear-Transport}), we will utilize the contraction mapping principle to prove the local well-posedness Theorem \ref{the-1}.
  To prove Theorem \ref{the-2},  we will show a  Kato-type estimate on the velocity and its gradient, which is
\begin{eqnarray*}
\| u(t)\|_{W^{1,\infty}(\R^2_+)}\leq C_\gamma\|\omega_0\|_{L^\infty(\R^2_+)}\Big(1+\log\Big(1+\frac{[\omega(t)]_{C^\gamma(\R^2_+)}}{\|\omega_0\|_{L^\infty(\R^2_+)}}\Big)\Big)
+C_p\|\omega_0\|_{L^p(\R^2_+)},
\end{eqnarray*}
for any $0<\gamma<1$ and $1<p<2$ (see Lemma \ref{Kato inequality}), and hence, for any $1<p<2$, by the Sobolev embedding $W^{2,p}(\R^2_+)\hookrightarrow C^{0,2-\frac{2}{p}}(\R^2_+)$, it holds that
\begin{eqnarray*}
\|u(t)\|_{W^{1,\infty}(\R^2_+)}\leq C_p(\|\omega_0\|_{L^\infty(\R^2_+)}+\|\omega_0\|_{L^p(\R^2_+)})\Big(1+\log\Big(e+\|\omega(t)\|_{W^{2,p}(\R^2_+)}\Big)\Big)
,
\end{eqnarray*}
which is not sufficient to obtain the uniform estimate on $\displaystyle\sup_{0\leq t<T}\|\omega(t)\|_{W^{k,p}(\R^2_+)}$ with $k\geq2$ and $1<p<2$ for any $T>0$ (see \eqref{5.14} and \eqref{5.16}).
This phenomenon is very different from the well-known global well-posedness to the Cauchy problem \eqref{1.1}-\eqref{1.2} or initial-boundary problem \eqref{1.1}-\eqref{1.3} in the Sobolev spaces, where the logarithmic-type estimate of the gradient of the velocity  is enough to prove the global regularity of the solutions to \eqref{1.1}-\eqref{1.2} or \eqref{1.1}-\eqref{1.3}.
It is easily seen that this difference between the velocity equations and the vorticity equation is caused by their different non-linear structures.
In particular, there is a new term $\|\nabla\omega(t)\|_{L^2(\R^2_+)}$ appearing in the energy estimate of $\|\omega(t)\|_{W^{k,p}(\R^2_+)}$ when utilizing the calculus inequalities to the vorticity equation.
To complete the proof of Theorem \ref{the-2}, the uniform estimate on $\int^t_0\|u(\tau)\|_{W^{1,\infty}(\R^2_+)}d\tau$ is then required.
To this end, we will derive the double exponential growth in time of the gradient of the vorticity in the half plane (see Proposition \ref{Double-Exponetial-Upper-Bound}), which in turn implies that, for any $t>0$,
\begin{eqnarray*}
\|u(t)\|_{W^{1,\infty}(\R^2_+)}\leq C_pABe^{C_pAt},
\end{eqnarray*}
where $A=\|\omega_0\|_{L^\infty(\R^2_+)}+\|\omega_0\|_{L^p(\R^2_+)}$ and $B=1+\log\Big(3+\frac{\|\nabla\omega_0\|_{L^\infty(\R^2_+)}}{\|\omega_0\|_{L^\infty(\R^2_+)}}\Big)$.
Moreover, to establish the Kato-type estimate, we will also show a Schauder-type estimate as
\begin{eqnarray*}
[\nabla u]_{C^\gamma(\bar{\R}^2_+)}\leq C_\gamma [\omega]_{C^\gamma(\bar{\R}^2_+)},
\end{eqnarray*}
for any $0<\gamma<1$ (see Proposition \ref{Schauder-type-estimate}).
The similar   estimate  in the whole plane can be directly derived by making use of  the expression of the gradient of the velocity  and the cancellation  of the singular kernel (see for instance Lemma 4.6 in \cite{[MB]}). However, in the half-plane case,    the cancellation property of the singular kernel does not hold  due to the presence of the  boundary and the proof of in the whole plane can not be applied. We will overcome this difficulty with help of the stream function. In particular, the local maximum estimate of the stream function will be derived (see \eqref{3.17}).
 Finally, it should be  remarked that the optimal growth of the vorticity in the half plane remains open (see \cite{[KS]} for the unit disc and \cite{[X]} for symmetric smooth bounded domains).

This paper is organized as follows. In Section $2$, we first give some notations and definitions on Sobolev spaces and H\"{o}lder spaces, and then recall some useful facts utilized later  in this paper.
Section 3 is devoted to some delicate estimates between the velocity and the vorticity.
The proof of Theorems \ref{the-1} and \ref{the-2} is given in Sections $4$ and $5$, respectively.

\section{Preliminaries}
In this section, we will present some notations and basic facts needed later.
\subsection{Notations}
 Some notations are introduced as follows. Let $\Omega$ be an open set in $\R^N$.  For $s\geq1$, we denote $L^s(\Omega)$  the usual function space consisting of measurable functions on $\Omega$ which are s-integrable, of which  norm  is defined by
\begin{eqnarray*}
\|f\|_{L^s(\Omega)}=\Big(\int_{\Omega}|f(x)|^s dx\Big)^{\frac{1}{s}}.
\end{eqnarray*}
When $f$ is a vector or matrix function the same notation will be used and the notation $|f|$ denotes the usual Euclidean norm. For $s=\infty$, $L^\infty(\Omega)$ denotes the Banach space of essentially bounded functions on $\Omega$, of which norm is
\begin{eqnarray*}
\|f\|_{L^\infty(\Omega)}=\inf\{c:|f(x)|\leq c {\rm~a.e. ~on}~\Omega \}.
\end{eqnarray*}
An useful characterization of the $L^s$ norm of a function is by its distribution function. More precisely, for any $s>0$ and $|f|^s\in L^1(\Omega)$, it holds that
\begin{eqnarray}\label{2.1}
\int_{\Omega}|f(x)|^sdx=s\int^\infty_0\tau^{s-1}\mu_f(\tau)d\tau,
\end{eqnarray}
where $\mu_f$ is the distribution function of $f$ defined by
\begin{eqnarray}\label{2.2}
\mu_f(\tau)=|\{x\in\Omega:|f(x)|>\tau\}|
\end{eqnarray}
for any $\tau>0$ (see, e.g., \cite{[GT]}).

For $m\in\mathbb{N}_+$ a positive integer  and $1\leq s\leq\infty$, we define the usual Sobolev spaces $W^{m,s}(\Omega)
=\{f\in L^s(\Omega): \partial^\alpha f\in L^s(\Omega), |\alpha|\leq m\}$, where $\partial^\alpha f$ is the weak (or distributional)   derivative. These spaces are equipped with the following norms:
\begin{eqnarray*}
\|f\|_{W^{m,s}(\Omega)}=\Big(\sum_{|\alpha|\leq m}\|\partial^\alpha u\|^s_{L^s(\Omega)}\Big)^{\frac{1}{s}}, 1\leq s<\infty,
\end{eqnarray*}
and
\begin{eqnarray*}
\|f\|_{W^{m,\infty}(\Omega)}=\sum_{|\alpha|\leq m}\|\partial^\alpha f\|_{L^\infty(\Omega)}.
\end{eqnarray*}
There will be no notational distinction between Sobolev spaces of scalar-valued and vector-valued functions.
One can consult \cite{[AF],[B],[E],[GT],[St]} and the references therein for more materials of Sobolev spaces. Here we focus on the half-space case $\Omega=\R^N_+$.
Some Sobolev embeddings are listed here for convenient applications. More precisely,
$W^{m,s}(\R^N_+)\hookrightarrow L^{\frac{Ns}{N-ms}}(\R^N_+)$ if $ms<N$ and $W^{m,s}(\R^N_+)\hookrightarrow C^k(\bar{\R}^N_+)$ if $ms>N$ (see, e.g., \cite{[B]}).

  For $0<\gamma<1$, a function  $f$ is uniformly H\"{o}lder continuous with exponent $\gamma$ in $\Omega$ if the quantity
\begin{eqnarray*}
[f]_{C^\gamma(\Omega)}=\sup_{x\neq y\in\Omega}\frac{|f(x)-f(y)|}{|x-y|^\gamma}
\end{eqnarray*}
is finite, and locally  H\"{o}lder continuous with exponent $\gamma$ in $\Omega$ if $f$ is uniformly H\"{o}lder continuous with exponent $\gamma$ on compact subset of $\Omega$. For $k$ a non-negative integer and $\gamma\in(0,1)$, the H\"{o}lder spaces $C^{k,\gamma}(\bar{\Omega})$~$(C^{k,\gamma}(\Omega))$ are defined as the subspaces of $C^k(\bar{\Omega})~(C^k(\Omega))$ consisting of functions whose $k-{\rm th}$ order partial derivatives are uniformly H\"{o}lder continuous (locally H\"{o}lder continuous) with exponent $\gamma$ in $\Omega$. Set
\begin{eqnarray*}
[f]_{C^k(\Omega)}=\|D^k f\|_{C^0(\Omega)}=\sup_{|\beta|=k}\sup_{\Omega}|D^\beta f|,~k=0,1,2...,
\end{eqnarray*}
\begin{eqnarray*}
[f]_{C^{k,\gamma}(\Omega)}=[D^k f]_{C^\gamma(\Omega)}=\sup_{|\beta|=k}[D^\beta f]_{C^\gamma(\Omega)},
\end{eqnarray*}
which are semi-norms in $C^k(\Omega)$ and $C^{k,\gamma}(\Omega)$ respectively. With these semi-norms, we can define the related norms
\begin{eqnarray*}
\|f\|_{C^k(\bar{\Omega})}=\sum^k_{j=0}[f]_{C^j(\Omega)}=\sum^k_{j=0}\|D^j f\|_{C^0(\Omega)},
\end{eqnarray*}
\begin{eqnarray*}
\|f\|_{C^{k,\gamma}(\bar{\Omega})}
&=&\|f\|_{C^k(\bar{\Omega})}+[f]_{C^{k,\gamma}(\Omega)}\\
&=&\|f\|_{C^k(\bar{\Omega})}+[D^k f]_{C^\gamma(\Omega)},
\end{eqnarray*}
on the spaces $C^{k}(\bar{\Omega})$, $C^{k,\gamma}(\bar{\Omega})$, respectively.

We also introduce the non-dimensional norms on $C^{k}(\bar{\Omega})$, $C^{k,\gamma}(\bar{\Omega})$. If $\Omega$ is bounded, with $d=\rm diam~\Omega$ (the diameter of $\Omega$), we set
\begin{eqnarray*}
\|f\|^\prime_{C^k(\bar{\Omega})}=\sum^k_{j=0}d^j[f]_{C^j(\Omega)}=\sum^k_{j=0}d^j\|D^j f\|_{C^0(\Omega)},
\end{eqnarray*}
\begin{eqnarray*}
\|f\|^\prime_{C^{k,\gamma}(\bar{\Omega})}
&=&\|f\|^\prime_{C^k(\bar{\Omega})}+d^{k+\gamma}[f]_{C^{k,\gamma}(\Omega)}\\
&=&\|f\|^\prime_{C^k(\bar{\Omega})}+d^{k+\gamma}[D^k f]_{C^\gamma(\Omega)}.
\end{eqnarray*}
 The spaces $C^{k}(\bar{\Omega})$, $C^{k,\gamma}(\bar{\Omega})$ equipped with the respective norms are Banach spaces. We refer readers to \cite{[GT]} for more details.

Throughout this paper, we will use $C$ to denote a generic positive constant, whose value may
change from line to line, and write $C({\alpha})$ or $C_{\alpha}$ to emphasize the dependence of a
constant on $\alpha$.
\subsection{Some basic facts}
In the following, we present some basic useful facts needed later.
We begin with the contraction mapping principle \cite{[GT]}, which will be used to prove the local existence of a solution to \eqref{1.4}.
\begin{lemma}\label{contration}
 Assume that $\mathcal{F}$ is a closed nonempty subset of Banach space $\mathcal{B}$ equipped with a norm $\|\cdot\|$ and that a mapping $\mathcal{T}:\mathcal{F}\rightarrow \mathcal{F}$ is contractive, that is,
\begin{eqnarray*}
\|\mathcal{T}x-\mathcal{T}y\|\leq\kappa\|x-y\|,
\end{eqnarray*}
for all $x,y\in\mathcal{F}$ and some $\kappa\in[0,1)$.
Then there exists a unique solution $x\in\mathcal{F}$ of the equation $\mathcal{T}x=x$.
\end{lemma}
Next we recall the  Hardy-Littlewood-Sobolev inequality \cite{[BCD]} needed later.
\begin{lemma}\label{Hardy}
Let $\beta\in(0,N)$ and $s,r\in(1,\infty)$ satisfy
\begin{eqnarray*}
\frac{1}{s}+\frac{\beta}{N}=1+\frac{1}{r}.
\end{eqnarray*}
A constant $C_{r,s}$ then exists such that
\begin{eqnarray*}
\||\cdot|^{-\beta}\ast f\|_{L^r(\R^N)}\leq C_{r,s}\|f\|_{L^s(\R^N)}.
\end{eqnarray*}
\end{lemma}
We continue with an integral inequality \cite{[MP]}, which will play a similar role as the Gronwall's inequality.
\begin{lemma}\label{Integral-Inequality}
Let $f\in C([0,T];\R_+)$ and $\varphi\in C(\R_+;\R_+)$ be a non-decreasing function, such that
\begin{eqnarray*}
f(t)\leq f(0)+\int^{t}_{0}\varphi(f(\tau))d\tau, ~t\leq T.
\end{eqnarray*}
Let $g=g(t)$ be a solution of the initial value problem
\begin{equation*}
\left\{\ba
&\frac{d}{dt}g(t)=\varphi(g(t)),\\
&g(0)=f(0). \ea\ \right.
\end{equation*}
Then
\begin{eqnarray*}
f(t)\leq g(t),  t\in[0,T].
\end{eqnarray*}
\end{lemma}
The next lemma is from the standard theory of elliptic partial differential equations of second order \cite{[GT]}. Let $\R^N_+$ denote the half-space, $x_N>0$, and $\mathbb{T}$ the hyperplane $x_N=0$;
$B_2=B_{2R}(x_0)$, $B_1=B_{R}(x_0)$ will be balls with center $x_0\in \bar{\R}^N_+$ and we let $B^+_2=B_2\cap\R^N_+$, $B^+_1=B_1\cap\R^N_+$.
\begin{lemma}\label{Schauder at the boundary}
Let $f\in C^2(B^+_2)\cap C^0(\bar{B}^+_2)$, $g\in C^\gamma(\bar{B}^+_2)$ with $0<\gamma<1$, satisfy $\Delta f=g$ in $B^+_2$, $f=0$ on $\mathbb{T}$. Then $f\in C^{2,\gamma}(\bar{B}^+_1)$ and we have
\begin{eqnarray*}
\|f\|^\prime_{C^{2,\gamma}(\bar{B}^+_1)}\leq C_{N,\gamma}(\|f\|_{C^{0}(\bar{B}^+_2)}+R^2\|g\|^\prime_{C^{0,\gamma}(\bar{B}^+_2)}).
\end{eqnarray*}
\end{lemma}
We conclude this section with the general calculus inequalities in the Sobolev spaces in the whole space and in the half-space. These inequalities are motivated by the well-known Kato-Ponce inequality which was originated in \cite{[Kato-Ponce]} and further generalized in \cite{[Kato3],[KPV],[L-D]}.
\begin{proposition}\label{Calculus-whole-space}
Suppose that $s\in(1,\infty)$ and $\alpha=(\alpha_1,\cdots,\alpha_N)\in\mathbb{N}^{N}$ with $|\alpha|=\displaystyle\sum_{j=1}^{N}\alpha_j>0$. Let $f,g$ be in $\mathcal{S}(\R^N)$, the Schwartz class. Then there exist constants $C's$ depending only on
$N, \alpha, s, s_1, s_2, s_3$ and $s_4$ such that
\begin{eqnarray}\label{2.3}
\|\partial^{\alpha}(fg)\|_{L^s(\R^N)}
\leq
C(\|D^{|\alpha|}f\|_{L^{s_1}(\R^N)}
\|g\|_{L^{s_2}(\R^N)}+\|f\|_{L^{s_4}(\R^N)}\|D^{|\alpha|}g\|_{L^{s_3}(\R^N)}),
\end{eqnarray}
and
\begin{eqnarray}\begin{split}\label{2.4}
&\|\partial^{\alpha}(fg)-f\partial^{\alpha}g\|_{L^s(\R^N)}\\
\leq& C\Big(\|D^{|\alpha|}f\|_{L^{s_1}(\R^N)}\|g\|_{L^{s_2}(\R^N)}+\|\nabla f\|_{L^{s_4}(\R^N)}
\|D^{|\alpha|-1}g\|_{L^{s_3}(\R^N)}\Big),
\end{split}
\end{eqnarray}
with $s_1, s_3\in(1,\infty)$ such that
\begin{eqnarray*}
\frac{1}{s}=\frac{1}{s_1}+\frac{1}{s_2}=\frac{1}{s_3}+\frac{1}{s_4}.
\end{eqnarray*}
\end{proposition}
\begin{remark}
The special cases of \eqref{2.3}-\eqref{2.4} with $s=2, s_1=s_3=2, s_2=s_4=\infty$  and of \eqref{2.4} with $s\in(1,\infty), s_1=s_3=s, s_2=s_4=\infty$  are proved in \cite{[MB]} and \cite{[L-D]}, respectively.
\end{remark}
\textbf{Proof of Proposition \ref{Calculus-whole-space}.}
When $|\alpha|=1$, \eqref{2.3} is a direct consequence of the Leibniz differential formula and the H\"{o}lder inequality.
For $|\alpha|\geq2$, the Leibniz differential formula gives
\begin{eqnarray}\label{2.5}
\partial^\alpha(fg)
=(\partial^\alpha f) g+\sum^{|\alpha|-1}_{|\beta|=1}C_{\alpha,\beta}\partial^\beta f\partial^{\alpha-\beta}g+f\partial^\alpha g.
\end{eqnarray}
For $1\leq|\beta|\leq|\alpha|-1$, by using the Gagliardo-Nirenberg inequality (see \cite{[N]}), we get
\begin{eqnarray}\label{2.6}
\|\partial^\beta f\|_{L^{p_1}(\R^N)}
\leq
C_{N,\alpha,\beta,s_1,s_4}\|D^{|\alpha|} f\|^{\frac{|\beta|}{|\alpha|}}_{L^{s_1}(\R^N)}\|f\|^{1-\frac{|\beta|}{|\alpha|}}_{L^{s_4}(\R^N)},
\end{eqnarray}
and
\begin{eqnarray}\label{2.7}
\|\partial^{\alpha-\beta}g\|_{L^{p_2}(\R^N)}
\leq
C_{N,\alpha,\beta,s_2,s_3}\|D^{|\alpha|}g\|^{1-\frac{|\beta|}{|\alpha|}}_{L^{s_3}(\R^N)}\|g\|^{\frac{|\beta|}{|\alpha|}}_{L^{s_2}(\R^N)},
\end{eqnarray}
where $p_1$ and $p_2$ are determined by
\begin{eqnarray*}
\frac{1}{p_1}=\frac{|\beta|}{|\alpha|}\frac{1}{s_1}+\Big(1-\frac{|\beta|}{|\alpha|}\Big)\frac{1}{s_4},
\end{eqnarray*}
and
\begin{eqnarray*}
\frac{1}{p_2}=\Big(1-\frac{|\beta|}{|\alpha|}\Big)\frac{1}{s_3}+\frac{|\beta|}{|\alpha|}\frac{1}{s_2}.
\end{eqnarray*}
Then it is easily seen that $\frac{1}{s}=\frac{1}{p_1}+\frac{1}{p_2}$. Therefore, utilizing the H\"{o}lder inequality and the Young inequality together with \eqref{2.5}-\eqref{2.7} yields that, for any $|\alpha|\geq2$,
\begin{eqnarray*}
\begin{split}
&\|\partial^{\alpha}(fg)\|_{L^s(\R^N)}\\
\leq&\|D^{|\alpha|}f\|_{L^{s_1}(\R^N)}
\|g\|_{L^{s_2}(\R^N)}+\sum^{|\alpha|-1}_{|\beta|=1}C_{\alpha,\beta}\|\partial^\beta f\|_{L^{p_1}(\R^N)}\|\partial^{\alpha-\beta}g\|_{L^{p_2}(\R^N)}\\
&\
+\|f\|_{L^{s_4}(\R^N)}\|D^{|\alpha|}g\|_{L^{s_3}(\R^N)}\\
\leq&\|D^{|\alpha|}f\|_{L^{s_1}(\R^N)}
\|g\|_{L^{s_2}(\R^N)}+\|f\|_{L^{s_4}(\R^N)}\|D^{|\alpha|}g\|_{L^{s_3}(\R^N)}\\
&\
+\sum^{|\alpha|-1}_{|\beta|=1}C_{N,\alpha,\beta,s_i}\Big(\|D^{|\alpha|} f\|_{L^{s_1}(\R^N)}
\|g\|_{L^{s_2}(\R^N)}\Big)^{\frac{|\beta|}{|\alpha|}}
\Big(\|f\|_{L^{s_4}(\R^N)}\|D^{|\alpha|}g\|_{L^{s_3}(\R^N)}\Big)^{1-\frac{|\beta|}{|\alpha|}}\\
\leq& C_{N,\alpha,\beta,s_i}(\|D^{|\alpha|}f\|_{L^{s_1}(\R^N)}
\|g\|_{L^{s_2}(\R^N)}+\|f\|_{L^{s_4}(\R^N)}\|D^{|\alpha|}g\|_{L^{s_3}(\R^N)}),
\end{split}
\end{eqnarray*}
which is exactly \eqref{2.3}.

Next we turn to the proof of \eqref{2.4}. When $|\alpha|=1 ~{\rm or} ~2$, \eqref{2.4}
is straightforward by using the Leibniz differential formula and the H\"{o}lder inequality. Thus we proceed to assume $|\alpha|\geq3$.
Similar to \eqref{2.5}, one writes
\begin{eqnarray}\label{2.8}
\partial^{\alpha}(fg)-f\partial^{\alpha}g
=(\partial^{\alpha}f)g
+\sum_{|\beta|=1}C_{\alpha,\beta}\partial^\beta f\partial^{\alpha-\beta}g
+\sum^{|\alpha|-1}_{|\beta|=2}C_{\alpha,\beta}\partial^\beta f\partial^{\alpha-\beta}g.
\end{eqnarray}
Furthermore, similar to \eqref{2.6}-\eqref{2.7}, we deduce that, for any $2\leq|\beta|\leq|\alpha|-1$,
\begin{eqnarray}\label{2.9}
\|\partial^\beta f\|_{L^{q_1}(\R^N)}
\leq
C_{N,\alpha,\beta,s_1,s_4}\|D^{|\alpha|} f\|^{\frac{|\beta|-1}{|\alpha|-1}}_{L^{s_1}(\R^N)}\|\nabla f\|^{1-\frac{|\beta|-1}{|\alpha|-1}}_{L^{s_4}(\R^N)},
\end{eqnarray}
and
\begin{eqnarray}\label{2.10}
\|\partial^{\alpha-\beta}g\|_{L^{q_2}(\R^N)}
\leq
C_{N,\alpha,\beta,s_2,s_3}\|D^{|\alpha|-1}g\|^{\frac{|\alpha|-|\beta|}{|\alpha|-1}}_{L^{s_3}(\R^N)}\|g\|^{1-\frac{|\alpha|-|\beta|}{|\alpha|-1}}_{L^{s_2}(\R^N)},
\end{eqnarray}
where $q_1$ and $q_2$ satisfy
\begin{eqnarray*}
\frac{1}{q_1}=\frac{|\beta|-1}{|\alpha|-1}\frac{1}{s_1}+\Big(1-\frac{|\beta|-1}{|\alpha|-1}\Big)\frac{1}{s_4},
\end{eqnarray*}
and
\begin{eqnarray*}
\frac{1}{q_2}=\frac{|\alpha|-|\beta|}{|\alpha|-1}\frac{1}{s_3}+\Big(1-\frac{|\alpha|-|\beta|}{|\alpha|-1}\Big)\frac{1}{s_2}.
\end{eqnarray*}
Then it is clear that $\frac{1}{s}=\frac{1}{q_1}+\frac{1}{q_2}$. Therefore, utilizing the H\"{o}lder inequality and Young inequality together with \eqref{2.8}-\eqref{2.10} deduces that, for any $|\alpha|\geq3$,
\begin{eqnarray*}
\begin{split}
&\|\partial^{\alpha}(fg)-f\partial^{\alpha}g\|_{L^s(\R^N)}\\
\leq&\|D^{|\alpha|}f\|_{L^{s_1}(\R^N)}
\|g\|_{L^{s_2}(\R^N)}+\|\nabla f\|_{L^{s_4}(\R^N)}\|D^{|\alpha|-1}g\|_{L^{s_3}(\R^N)}\\
&\
+\sum_{2\leq|\beta|\leq|\alpha|-1}C_{\alpha,\beta}\|\partial^\beta f\|_{L^{q_1}(\R^N)}\|\partial^{\alpha-\beta}g\|_{L^{q_2}(\R^N)}\\
\leq&\|D^{|\alpha|}f\|_{L^{s_1}(\R^N)}\|g\|_{L^{s_2}(\R^N)}+\|\nabla f\|_{L^{s_4}(\R^N)}\|D^{|\alpha|-1}g\|_{L^{s_3}(\R^N)}\\
&\
+\sum^{|\alpha|-1}_{|\beta|=2}C_{N,\alpha,\beta,s_i}(\|D^{|\alpha|} f\|_{L^{s_1}(\R^N)}
\|g\|_{L^{s_2}(\R^N)})^{\frac{|\beta|-1}{|\alpha|-1}}
(\|\nabla f\|_{L^{s_4}(\R^N)}\|D^{|\alpha|-1}g\|_{L^{s_3}(\R^N)})^{\frac{|\alpha|-|\beta|}{|\alpha|-1}}\\
\leq& C_{N,\alpha,\beta,s_i}(\|D^{|\alpha|}f\|_{L^{s_1}(\R^N)}
\|g\|_{L^{s_2}(\R^N)}+\|\nabla f\|_{L^{s_4}(\R^N)}\|D^{|\alpha|-1}g\|_{L^{s_3}(\R^N)}),
\end{split}
\end{eqnarray*}
which concludes the proof of \eqref{2.4}.

The corresponding calculus inequalities in the half-space are then easily derived from Stein's linear Sobolev extension operators (\cite{[St]}) and Proposition \ref{Calculus-whole-space}.
\begin{coro}\label{Calculus-half-space}
Suppose that $\alpha=(\alpha_1,\cdots,\alpha_N)\in\mathbb{N}^{N}$ with $|\alpha|=\displaystyle\sum_{j=1}^{N}\alpha_j>0$ and $s\in(1,\infty)$ with $s_1, s_3\in(1,\infty)$ such that $\frac{1}{s}=\frac{1}{s_1}+\frac{1}{s_2}=\frac{1}{s_3}+\frac{1}{s_4}$.

{\rm(1)} If $f\in W^{|\alpha|,s_1}(\R^N_+)\cap L^{s_4}(\R^N_+)$ and $g\in W^{|\alpha|,s_3}(\R^N_+)\cap L^{s_2}(\R^N_+)$, then there exists a constant $C=C(N, \alpha, s, s_1, s_2, s_3,s_4)$ such that
\begin{eqnarray*}
\|fg\|_{W^{|\alpha|,s}(\R^N_+)}
\leq C(\|f\|_{W^{|\alpha|,s_1}(\R^N_+)}\|g\|_{L^{s_2}(\R^N_+)}+\|f\|_{L^{s_4}(\R^N_+)}\|g\|_{W^{|\alpha|,s_3}(\R^N_+)}).
\end{eqnarray*}

{\rm(2)} If $f\in W^{|\alpha|,s_1}(\R^N_+)\cap W^{1,s_4}(\R^N_+)$ and $g\in W^{|\alpha|-1,s_3}(\R^N_+)\cap L^{s_2}(\R^N_+)$, then there exists a constant $C=C(N, \alpha, s, s_1, s_2, s_3,s_4)$ such that
\begin{eqnarray*}
\|\partial^{\alpha}(fg)-f\partial^{\alpha}g\|_{L^s(\R^N_+)}
\leq C(\|f\|_{W^{|\alpha|,s_1}(\R^N_+)}\|g\|_{L^{s_2}(\R^N_+)}+\|f\|_{W^{1,s_4}(\R^N_+)}
\|g\|_{W^{|\alpha|-1,s_3}(\R^N_+)}).
\end{eqnarray*}
\end{coro}
\section{Estimates on the velocity fields}
Let us start with some estimates on the  velocity field $u$  in the half plane through the vorticity $\omega$. These estimates are of interest in itself. For convenience, we will make a notation suppression for the time variable in this section. All functions could be understood at a fixed time $t>0$.

As usual,  the stream function $\Psi(x)$, which vanishes at infinity, is introduced and satisfies
\begin{equation}\label{3.1}
\left\{\ba
&\Delta\Psi(x)=\omega(x),~x\in\R^2_+,\\\
&\Psi(x_1,0)=0.\ea\ \right.
\end{equation}
 Clearly, the stream function $\Psi$ can be represented as
\begin{eqnarray}\label{3.2}
\Psi(x)=\frac{1}{2\pi}\int_{\mathbb R^2_+}\Big(\log|x-y|-\log|x-\bar{y}|\Big)\omega(y)dy,~x\in\R^2_+.
\end{eqnarray}
And the velocity field $u=\nabla^{\bot}\Psi$ can be explicitly formulated as
\begin{eqnarray}\label{3.3}
u(x)=\frac{1}{2\pi}\int_{\mathbb R^2_+}\Big(\frac{(x-y)^\bot}{|x-y|^2}-\frac{(x-\bar{y})^\bot}{|x-\bar{y}|^2}\Big)\omega(y)dy,~x\in\R^2_+,
\end{eqnarray}
which is called the Biot-Savart law. We use here the notations $z^\perp:=(z_2,-z_1)$ and $\bar{z}:=(z_1,-z_2)$ for $z=(z_1,z_2).$
Then the velocity field $u$ satisfies
\begin{eqnarray}\label{3.4}
{\rm div}~u=0,~x\in\R^2_+~{\rm and}~ u_2(x_1,0)=0, ~x_1\in\R.
\end{eqnarray}
Indeed,  the slip boundary condition  and the divergence-free condition in \eqref{3.4} can be derived from \eqref{3.3} and the expression  of the gradient of the velocity, respectively. In particular, the expression  of the gradient of the velocity is calculated  in  Lemma \ref{Formula for gradient}.

The following estimate is a Sobolev-type estimate of the velocity field $u$ in terms of the vorticity $\omega$ by the Biot-Savart law.
\begin{proposition}\label{Sobolev-type estimate}
Let $\omega\in W^{m,p}(\R^2_+)$ with  $m\geq0$ an integer and $1<p<2$. Suppose that $u$ is defined by \eqref{3.3}. Then
\begin{eqnarray*}
\|u\|_{W^{m,\frac{2p}{2-p}}(\R^2_+)}\leq C_{m,p}\|\omega\|_{W^{m,p}(\R^2_+)}.
\end{eqnarray*}
\end{proposition}
{\bf Proof.}
When $m=0$, from \eqref{3.3}, we have
\begin{eqnarray*}
u(x)
&=&\frac{1}{2\pi}\int_{\mathbb R^2_+}\Big(\frac{(x-y)^\bot}{|x-y|^2}-\frac{(x-\bar{y})^\bot}{|x-\bar{y}|^2}\Big)\omega(y)dy\\
&=&\frac{1}{2\pi}\int_{\mathbb R^2}\frac{(x-y)^\bot}{|x-y|^2}\bar{\omega}(y)dy,
\end{eqnarray*}
where $\bar{\omega}:\R^2\rightarrow \R $ is the odd extension of $\omega$ to the whole plane, which is
\begin{eqnarray*}
\bar{\omega}(y_1,y_2)=
\begin{cases}
\omega(y_1,y_2),
&  y_2>0, \\
-\omega(y_1,-y_2),
& y_2<0. \\
\end{cases}
\end{eqnarray*}
Hence, for any $x\in \R^2_+$,
\begin{eqnarray*}
|u(x)|
&\leq&\frac{1}{2\pi}\int_{\mathbb R^2}\frac{1}{|x-y|}|\bar{\omega}(y)|dy\\
&=&\frac{1}{2\pi}(|\cdot|^{-1}\ast|\bar{\omega}|)(x).
\end{eqnarray*}
It follows from Lemma \ref{Hardy} that, for $p\in(1,2)$,
\begin{eqnarray}\label{3.5}
\|u\|_{L^{\frac{2p}{2-p}}(\R^2_+)}
&\leq&\|\frac{1}{2\pi}|\cdot|^{-1}\ast|\bar{\omega}|\|_{L^{\frac{2p}{2-p}}(\R^2_+)}\nonumber\\
&\leq&\|\frac{1}{2\pi}|\cdot|^{-1}\ast|\bar{\omega}|\|_{L^{\frac{2p}{2-p}}(\R^2)}\nonumber\\
&\leq&C_p\|\bar{\omega}\|_{L^p(\R^2)}\nonumber\\
&\leq&C_p\|\omega\|_{L^p(\R^2_+)}.
\end{eqnarray}
We proceed to show the case $m=1$.
From Lemma \ref{Formula for gradient}, we see that
\begin{eqnarray*}
\nabla u(x)=\frac{1}{2\pi}P.V.\int_{\R^2_+}M(x,y)\omega(y)dy+\frac{\omega(x)}{2}
\begin{pmatrix}
0&1\\-1&0
\end{pmatrix},
\end{eqnarray*}
where the kernel matrix $M(x,y)$ is the sum of singular  part $M_s(x,y)$ and regular part $M_r(x,y)$:
\begin{eqnarray*}
M_s(x,y)=
\begin{pmatrix}
\frac{-2(x_1-y_1)(x_2-y_2)}{|x-y|^4}
&\frac{(x_1-y_1)^2-(x_2-y_2)^2}{|x-y|^4}\\
\frac{(x_1-y_1)^2-(x_2-y_2)^2}{|x-y|^4}
&\frac{2(x_1-y_1)(x_2-y_2)}{|x-y|^4}
\end{pmatrix},
\end{eqnarray*}
\begin{eqnarray*}
M_r(x,y)=
\begin{pmatrix}
\frac{2(x_1-y_1)(x_2+y_2)}{|x-\bar{y}|^4}
&-\frac{(x_1-y_1)^2-(x_2+y_2)^2}{|x-\bar{y}|^4}\\
-\frac{(x_1-y_1)^2-(x_2+y_2)^2}{|x-\bar{y}|^4}
&-\frac{2(x_1-y_1)(x_2+y_2)}{|x-\bar{y}|^4}
\end{pmatrix}.
\end{eqnarray*}
Utilizing the same odd extension technique, we can obtain
\begin{eqnarray*}
\nabla u(x)=\frac{1}{2\pi}P.V.\int_{\R^2}M_s(x,y)\bar{\omega}(y)dy+\frac{\omega(x)}{2}
\begin{pmatrix}
0&1\\-1&0
\end{pmatrix}.
\end{eqnarray*}
Therefore, by the classical Calder\'{o}n-Zygmund singular integral theory (see, e.g., \cite{[St]}), we have
\begin{eqnarray}\label{3.6}
\|\nabla u\|_{L^{\frac{2p}{2-p}}(\R^2_+)}
&\leq& C_p\|\bar{\omega}\|_{L^{\frac{2p}{2-p}}(\R^2)}+C\|\omega\|_{L^{\frac{2p}{2-p}}(\R^2_+)}\nonumber\\
&\leq& C_p\|\omega\|_{L^{\frac{2p}{2-p}}(\R^2_+)}\nonumber\\
&\leq& C_p\|\omega\|_{W^{1,p}(\R^2_+)},
\end{eqnarray}
where we have used the Sobolev embedding $W^{1,p}(\R^2_+)\hookrightarrow L^{\frac{2p}{2-p}}(\R^2_+)$ for any $p\in(1,2)$.
Combining \eqref{3.5} with \eqref{3.6} leads to
\begin{eqnarray}\label{3.7}
\|u\|_{W^{1,\frac{2p}{2-p}}(\R^2_+)}\leq C_p\|\omega\|_{W^{1,p}(\R^2_+)}.
\end{eqnarray}
We now turn to the case $m=2$. In view of \eqref{3.7}, it suffices to bound $\|\nabla^2u\|_{L^{\frac{2p}{2-p}}(\R^2_+)}$. To this end,  we first estimate the tangential derivative $\partial_{x_1}\nabla u$ of the gradient velocity $\nabla u$, that is, $\nabla_x\nabla^\bot_x\partial_{x_1}\Psi$.
After differentiating \eqref{3.2} directly and integrating by parts, together with the homogeneous boundary condition of the Green's function, we can obtain
\begin{eqnarray}\label{3.8}
\partial_{x_1}\Psi(x,t)
&=&\frac{1}{2\pi}\int_{\mathbb R^2_+}\partial_{x_1}\Big(\log|x-y|-\log|x-\bar{y}|\Big)\omega(y)dy\nonumber\\
&=&-\frac{1}{2\pi}\int_{\mathbb R^2_+}\partial_{y_1}\Big(\log|x-y|-\log|x-\bar{y}|\Big)\omega(y)dy\nonumber\\
&=&\frac{1}{2\pi}\int_{\mathbb R^2_+}\Big(\log|x-y|-\log|x-\bar{y}|\Big)\partial_{y_1}\omega(y)dy.
\end{eqnarray}
By differentiating \eqref{3.8}, we can arrive at
\begin{eqnarray*}
\nabla^\bot_x\partial_{x_1}\Psi(x)=\frac{1}{2\pi}\int_{\mathbb R^2_+}\Big(\frac{(x-y)^\bot}{|x-y|^2}-\frac{(x-\bar{y})^\bot}{|x-\bar{y}|^2}\Big)\partial_{y_1}\omega(y)dy.
\end{eqnarray*}
Computations similar to those in the proof of Lemma \ref{Formula for gradient} give
\begin{eqnarray*}
\nabla_x\nabla^\bot_x\partial_{x_1}\Psi(x)=\frac{1}{2\pi}P.V.\int_{\R^2_+}M(x,y)\partial_{y_1}\omega(y)dy+\frac{\partial_{x_1}\omega(x)}{2}
\begin{pmatrix}
0&1\\-1&0
\end{pmatrix}.
\end{eqnarray*}
Similar to \eqref{3.6}, we can obtain
\begin{eqnarray}\label{3.9}
\|\partial_{x_1}\nabla u\|_{L^{\frac{2p}{2-p}}(\R^2_+)}
&=&\|\nabla_x\nabla^\bot_x\partial_{x_1}\Psi\|_{L^{\frac{2p}{2-p}}(\R^2_+)}\nonumber\\
&\leq& C_p\|\partial_{x_1}\omega\|_{L^{\frac{2p}{2-p}}(\R^2_+)}\nonumber\\
&\leq& C_p\|\omega\|_{W^{2,p}(\R^2_+)}.
\end{eqnarray}
Next we consider the normal derivative $\partial_{x_2}\nabla u$, that is, $\partial_{x_2}\nabla_x\nabla_x^\bot\Psi$. In view of \eqref{3.9}, it remains to bound $\|\partial^3_{x_2}\Psi\|_{L^{\frac{2p}{2-p}}(\R^2_+)}$.
By \eqref{3.1}, we proceed to represent $\partial^3_{x_2}\Psi$ as
\begin{eqnarray*}
\partial^3_{x_2}\Psi=\partial_{x_2}(\omega-\partial^2_{x_1}\Psi).
\end{eqnarray*}
Using \eqref{3.9} again, we deduce that
\begin{eqnarray}\label{3.10}
\|\partial^3_{x_2}\Psi\|_{L^{\frac{2p}{2-p}}(\R^2_+)}
&\leq&\|\partial_{x_2}\omega\|_{L^{\frac{2p}{2-p}}(\R^2_+)}
+\|\partial^2_{x_1}\partial_{x_2}\Psi\|_{L^{\frac{2p}{2-p}}(\R^2_+)}\nonumber\\
&\leq& C_p\|\omega\|_{W^{2,p}(\R^2_+)}.
\end{eqnarray}
Finally, \eqref{3.9} and \eqref{3.10} yield
\begin{eqnarray*}
\|\nabla^2u\|_{L^{\frac{2p}{2-p}}(\R^2_+)}\leq C_p\|\omega\|_{W^{2,p}(\R^2_+)},
\end{eqnarray*}
which is the desired conclusion.
The case $m\geq3$ can be proved by the induction argument. We omit the details for simplicity and the proof of Proposition \ref{Sobolev-type estimate} is complete.

The following estimate is a Schauder-type one on the gradient of the velocity, which will be used to derive its maximum bound  in Lemma \ref{Kato inequality}. The estimate in the whole plane can be found in \cite{[MB]}.
\begin{proposition}\label{Schauder-type-estimate}
Let $\omega\in C^\gamma(\bar{\R}^2_+)\cap L^q(\R^2_+)$  with some $0<\gamma<1$ and $1\leq q<2$. Suppose that $\Psi$ is defined by \eqref{3.2} and $u=\nabla^{\bot}\Psi$. Then the Hessian $\nabla^2\Psi$ satisfies
\begin{eqnarray*}
[\nabla^2\Psi]_{C^\gamma(\bar{\R}^2_+)}\leq C_\gamma [\omega]_{C^\gamma(\bar{\R}^2_+)}.
\end{eqnarray*}
Consequently,
\begin{eqnarray*}
[\nabla u]_{C^\gamma(\bar{\R}^2_+)}\leq C_\gamma [\omega]_{C^\gamma(\bar{\R}^2_+)}.
\end{eqnarray*}
\end{proposition}
\begin{remark}
The condition $\omega\in L^q(\R^2_+)$ with $1\leq q<2$ in Proposition \ref{Schauder-type-estimate} is needed to obtain the local maximum estimate of the stream function $\Psi$ {\rm(}see \eqref{3.17}{\rm)}. But the H\"{o}lder semi-norm estimates of the Hessian $\nabla^2\Psi$ and the gradient velocity $\nabla u$ do not depend on the quantity $\|\omega\|_{L^q(\R^2_+)}$. Therefore the integrability condition $\omega\in L^q(\R^2_+)$ with $1\leq q<2$ may be removed by a suitable density argument.
\end{remark}

{\bf Proof of Proposition \ref{Schauder-type-estimate}.} Since
$u=\nabla^{\bot}\Psi$, we have $[\nabla u]_{C^\gamma(\bar{\R}^2_+)}\leq [\nabla^2\Psi]_{C^\gamma(\bar{\R}^2_+)}$.
Thus it suffices to estimate $[\nabla^2\Psi]_{C^\gamma(\bar{\R}^2_+)}$.
For any $R>0$,
it follows from \eqref{3.1} that
\begin{equation*}
\left\{\ba
&\Delta\Psi(x)=\omega(x),~x=(x_1,x_2)\in B^+_{2R}(0),\\\
&\Psi(x_1,0)=0.\ea\ \right.
\end{equation*}
By the definition of $\|\Psi\|^\prime_{C^{2,\gamma}(\bar{B}^+_{R}(0))}$ and Lemma \ref{Schauder at the boundary}, we obtain
\begin{eqnarray}\label{3.11}
&&(2R)^{2+\gamma}[\nabla^2\Psi]_{C^\gamma(\bar{B}^+_{R}(0))}\nonumber\\
&\leq&\|\Psi\|^\prime_{C^{2,\gamma}(\bar{B}^+_{R}(0))}\nonumber\\
&\leq&C_{\gamma}\Big(\|\Psi\|_{C^0(\bar{B}^+_{2R}(0))}+R^2\|\omega\|^\prime_{C^{0,\gamma}(\bar{B}^+_{2R}(0))}\Big)\nonumber\\
&=&C_{\gamma}\Big(\|\Psi\|_{C^0(\bar{B}^+_{2R}(0))}+R^2\|\omega\|_{C^0(\bar{B}^+_{2R}(0))}+R^{2+\gamma}[\omega]_{C^\gamma(\bar{B}^+_{2R}(0))}\Big)\nonumber\\
&\leq&C_{\gamma}\Big(\|\Psi\|_{C^0(\bar{B}^+_{2R}(0))}+R^2\|\omega\|_{C^0(\bar{\R}^2_+)}+R^{2+\gamma}[\omega]_{C^\gamma(\bar{\R}^2_+)}\Big).
\end{eqnarray}
Dividing the factor $(2R)^{2+\gamma}$ on  both sides of \eqref{3.11} implies
\begin{eqnarray}\label{3.12}
[\nabla^2\Psi]_{C^\gamma(\bar{B}^+_{R}(0))}\leq C_{\gamma}\Big(\frac{1}{R^{2+\gamma}}\|\Psi\|_{C^0(\bar{B}^+_{2R}(0))}+\frac{1}{R^\gamma}\|\omega\|_{C^0(\bar{\R}^2_+)}
+[\omega]_{C^\gamma(\bar{\R}^2_+)}\Big).
\end{eqnarray}
We proceed to estimate the term $\|\Psi\|_{C^0(\bar{B}^+_{2R}(0))}$ by utilizing the explicit expression \eqref{3.2}.
For any $x=(x_1,x_2)\in \bar{B}^+_{2R}(0)$, one easily has that $x_2\leq2R$. We rewrite the expression \eqref{3.2} as
\begin{eqnarray}\label{3.13}
\Psi(x)
&=&\frac{1}{2\pi}\int_{\R^2_+}\Big(\log|x-y|-\log|x-\bar{y}|\Big)\omega(y)dy\nonumber\\
&=&\frac{1}{2\pi}\Big(\int_{\R^2_+\cap B_1(x)}+\int_{\R^2_+\cap (B_1(x))^c}\Big)\Big(\log|x-y|-\log|x-\bar{y}|\Big)\omega(y)dy\nonumber\\
&=&\frac{1}{2\pi}\int_{\R^2_+\cap B_1(x)}\log|x-y|\omega(y)dy-\frac{1}{2\pi}\int_{\R^2_+\cap B_1(x)}\log|x-\bar{y}|\omega(y)dy\nonumber\\
&&\
+\frac{1}{2\pi}\int_{\R^2_+\cap (B_1(x))^c}\Big(\log|x-y|-\log|x-\bar{y}|\Big)\omega(y)dy\nonumber\\
&:=&J_{11}+J_{12}+J_{2}.
\end{eqnarray}
 $J_{11}$ is estimated as
\begin{eqnarray}\label{3.14}
|J_{11}|
&\leq&\frac{1}{2\pi}\int_{\R^2_+\cap B_1(x)}|\log|x-y|||\omega(y)|dy\nonumber\\
&\leq&-\frac{1}{2\pi}\|\omega\|_{L^\infty(\R^2_+)}\int_{B_1(x)}\log|x-y|dy\nonumber\\
&=&-\|\omega\|_{L^\infty(\R^2_+)}\int^1_0r\log r dr\nonumber\\
&=&C\|\omega\|_{L^\infty(\R^2_+)}.
\end{eqnarray}
$J_{12}$ is estimated as
\begin{eqnarray}\label{3.15}
|J_{12}|
&\leq& \frac{1}{2\pi}\int_{\R^2_+\cap B_1(x)}|\log|x-\bar{y}|||\omega(y)|dy\nonumber\\
&=& \frac{1}{2\pi}\int_{\R^2_+\cap B_1(x)\cap B_1(\bar{x})}|\log|x-\bar{y}|||\omega(y)|dy\nonumber\\
&&\
+ \frac{1}{2\pi}\int_{\R^2_+\cap B_1(x)\cap(B_1(\bar{x}))^c}|\log|x-\bar{y}|||\omega(y)|dy\nonumber\\
&=& -\frac{1}{2\pi}\int_{\R^2_+\cap B_1(x)\cap B_1(\bar{x})}\log|\bar{x}-y||\omega(y)|dy\nonumber\\
&&\
+ \frac{1}{2\pi}\int_{\R^2_+\cap B_1(x)\cap(B_1(\bar{x}))^c}\log|x-\bar{y}||\omega(y)|dy\nonumber\\
&\leq& -\frac{1}{2\pi}\int_{\R^2_+\cap B_1(\bar{x})}\log|\bar{x}-y||\omega(y)|dy\nonumber\\
&&\
+ \frac{1}{2\pi}\int_{\R^2_+\cap B_1(x)\cap(B_1(\bar{x}))^c}\log\Big(|x-\bar{x}|+|\bar{x}-\bar{y}|\Big)|\omega(y)|dy\nonumber\\
&=& -\frac{1}{2\pi}\int_{\R^2_+\cap B_1(\bar{x})}\log|\bar{x}-y||\omega(y)|dy\nonumber\\
&&\
+ \frac{1}{2\pi}\int_{\R^2_+\cap B_1(x)\cap(B_1(\bar{x}))^c}\log\Big(2x_2+|x-y|\Big)|\omega(y)|dy\nonumber\\
&\leq& -\|\omega\|_{L^\infty(\R^2_+)}\int^1_{0}r\log r dr
+ \frac{1}{2}\|\omega\|_{L^\infty(\R^2_+)}\log(1+4R)\nonumber\\
&\leq&C\Big(1+\log(1+4R)\Big)\|\omega\|_{L^\infty(\R^2_+)}.
\end{eqnarray}
For $J_2$, utilizing $\log(1+t)\leq t$ for $t\geq0$ and H\"{o}lder inequality, we obtain
\begin{eqnarray}\label{3.16}
|J_2|
&\leq&\frac{1}{4\pi}\int_{\R^2_+\cap (B_1(x))^c}\log\frac{|x-\bar{y}|^2}{|x-y|^2}|\omega(y)|dy\nonumber\\
&=&\frac{1}{4\pi}\int_{\R^2_+\cap (B_1(x))^c}\log\Big(1+\frac{4x_2y_2}{|x-y|^2}\Big)|\omega(y)|dy\nonumber\\
&\leq&\frac{1}{\pi}\int_{\R^2_+\cap (B_1(x))^c}\frac{x_2y_2}{|x-y|^2}|\omega(y)|dy\nonumber\\
&=&\frac{1}{\pi}\int_{\R^2_+\cap (B_1(x))^c}\frac{x_2(y_2-x_2)+x^2_2}{|x-y|^2}|\omega(y)|dy\nonumber\\
&\leq&\frac{1}{\pi}\int_{\R^2_+\cap (B_1(x))^c}\Big(\frac{x_2}{|x-y|}+\frac{x^2_2}{|x-y|^2}\Big)|\omega(y)|dy\nonumber\\
&\leq&\frac{2}{\pi}(R+2R^2)\int_{\R^2_+\cap (B_1(x))^c}\frac{1}{|x-y|}|\omega(y)|dy\nonumber\\
&\leq&C_q(R+R^2)\|\omega\|_{L^q(\R^2_+)},
\end{eqnarray}
for any $1\leq q<2$.
Consequently, combining \eqref{3.13}-\eqref{3.16} leads to
\begin{eqnarray}\label{3.17}
\|\Psi\|_{C^0(\bar{B}^+_R(0))}\leq C_q(R+R^2)\|\omega\|_{L^q(\R^2_+)}+C(1+\log(1+4R))\|\omega\|_{L^\infty(\R^2_+)},
\end{eqnarray}
which together with \eqref{3.12} deduces that, for any $R>0$,
\begin{eqnarray}\label{3.18}
\left.
 \begin{array}{ll}
[\nabla^2\Psi]_{C^\gamma(\bar{B}^+_{R}(0))}
&\leq C_{\gamma}\Big(C_q\Big(\frac{1}{R^\gamma}+\frac{1}{R^{1+\gamma}}\Big)\|\omega\|_{L^q(\R^2_+)} \\[3mm]
&\ \ \ \displaystyle
+\Big(\frac{1}{R^\gamma}+\frac{1}{R^{2+\gamma}}(1+\log(1+4R))\Big)\|\omega\|_{C^0(\bar{\R}^2_+)}
+[\omega]_{C^\gamma(\bar{\R}^2_+)}\Big).
 \end{array}
\right.
\end{eqnarray}
Letting $R\rightarrow\infty$ in \eqref{3.18} yields
\begin{eqnarray*}
[\nabla^2\Psi]_{C^\gamma(\bar{\R}^2_+)}\leq C_\gamma[\omega]_{C^\gamma(\bar{\R}^2_+)},
\end{eqnarray*}
which concludes the proof of Proposition \ref{Schauder-type-estimate}.
\section{Proof of Theorem \ref{the-1}}
In this section, we will prove Theorem \ref{the-1} by using the contraction mapping principle.
Before that, we introduce a nonempty closed subset contained in $C([0,T_0];L^{p}(\R^2_+))$ for some $1<p<2$.
Given $\omega_{0}\in W^{k,p}(\R^{2}_+)$ with $k\geq3$ and $1<p<2$, we define
\begin{equation}\begin{split}\label{4.1}
B_{T_0}=&\{\omega\in L^\infty(0,T_0;W^{k,p}(\R^{2}_+))\cap C([0,T_0];L^{p}(\R^2_+)):
ess\sup_{0\leq t\leq T_0} \|\omega(t)\|_{W^{k,p}(\R^{2}_+)}\leq M,\\&
 ~\omega(\cdot,0)=\omega_{0}\in W^{k,p}(\R^{2}_+)\},
 \end{split}
\end{equation}
where $M=2\|\omega_{0}\|_{W^{k,p}(\R^{2}_+)}$ and $T_0>0$ is  to be determined later.

We will construct the solution to \eqref{1.4} as a fixed point of a mapping $\mathcal{T}:B_{T_0}\rightarrow B_{T_0}$.
For any $\omega\in B_{T_0}$,
let us consider a linear transport equation in the upper half plane
\begin{equation}\label{4.2}
\left\{\ba
&\partial_t\theta+u\cdot\nabla \theta=0,\\
&\theta(\cdot,0)=\omega_{0}, \ea\ \right.
\end{equation}
where the advective velocity field  $u$ is determined by \eqref{3.3}.
The following lemma is concerned with the global existence and uniqueness of solutions to \eqref{4.2},  of which proof will be given in the Appendix A.
\begin{lemma}\label{Linear-Transport}
Let $\omega_{0}\in W^{k,p}(\R^2_+)$ with $k\geq3$ and $1<p<2$. Then for every $T>0$ and $\omega\in L^\infty(0,T;W^{k,p}(\R^2_+))\cap C([0,T];L^p(\R^2_+))$ satisfying $\omega(\cdot,0)=\omega_{0}$,   there exists a unique solution $\theta\in L^\infty(0,T;W^{k,p}(\R^2_+))\cap{\rm Lip}([0,T];W^{k-1,p}(\R^2_+))$ to \eqref{4.2}.
Moreover, it holds that
\begin{eqnarray}\label{4.3}
\|\theta(t)\|_{W^{k,p}(\R^2_+)}
\leq\|\omega_0\|_{W^{k,p}(\R^2_+)}
e^{C{\int_{0}^{t}\|\omega(\tau)\|_{W^{k,p}(\R^2_+)}d\tau}}
\end{eqnarray}
 for a.e. $t\in[0,T]$.
\end{lemma}
According to Lemma \ref{Linear-Transport}, we can define a mapping
\begin{eqnarray}\label{4.4}
\mathcal{T}(\omega)(x,t)=\theta(x,t),
\end{eqnarray}
which is from $L^\infty(0,T;W^{k,p}(\R^2_+))\cap C([0,T];L^p(\R^2_+))$ to $\theta\in L^\infty(0,T;W^{k,p}(\R^2_+))\cap{\rm Lip}([0,T];W^{k-1,p}(\R^2_+))$ for any $T>0$. Then we can show that
\begin{lemma}\label{existence of fixed-point}
The mapping $\mathcal{T}$ defined by \eqref{4.4} maps $B_{T_0}$ into itself for some $T_0>0$ and has exactly one fixed point in $B_{T_0}$.
\end{lemma}
{\bf Proof.} The proof is divided into three steps.

{\it Step 1. $\mathcal{T}$ maps $B_{T_0}$ into $B_{T_0}$.}

For $\omega\in B_{T_0}$, in view of Lemma \ref{Linear-Transport}, it concludes that $
 \theta(x,t)=\mathcal{T}(\omega)(x,t)\in C([0,T_0];L^{p}(\R^2_+))$ satisfying $\theta(\cdot,0)=\omega_0\in W^{k,p}(\R^{2}_+).$
Besides, by \eqref{4.3}, we can obtain
\begin{eqnarray*}
\|\theta\|_{L^\infty(0,T_0;W^{k,p}(\R^2_+))}&\leq&\|\omega_0\|_{W^{k,p}(\R^2_+)}\exp\{CT_0\|\omega\|_{L^\infty(0,T_0; W^{k,p}(\R^2_+))}\}\\
&\leq&\frac{1}{2}Me^{CMT_0}\\
&\leq&M,
\end{eqnarray*}
provided $0<T_0\le \frac{\log2}{CM}$.
This shows that $\theta\in B_{T_0}.$\\[2mm]
{\it Step 2. $B_{T_0}$ is a closed nonempty subset of $C([0,T_0];L^{p}(\R^2_+))$.}

$B_{T_0}$ is nonempty since $\omega_0\in B_{T_0}$.
To show that $B_{T_0}$ is closed in the topology of $C([0,T_0];L^{p}(\R^2_+))$, we assume that $\omega_n \in B_{T_0}$ and
\begin{eqnarray*}
\|\omega_n-\omega\|_{C([0,T_0];L^{p}(\R^2_+))}\rightarrow0, ~{\rm as}~n\rightarrow\infty.
\end{eqnarray*}
We have to show that  $\omega\in B_{T_0}$.
Since $\omega_n$ is bounded in $L^\infty(0,T_0;W^{k,p}(\R^2_+))$, that is, $\|\omega_n\|_{L^\infty(0,T_0;W^{k,p}(\R^2_+))}\leq M$, by the weak-star compactness, we can extract a subsequence $\{\omega_{n_j}\}_{j\in \mathbb{N}}$ such that $\omega_{n_j}\overset{\ast}{\rightharpoonup}\widetilde{\omega}$ in $L^\infty(0,T_0;W^{k,p}(\R^2_+)).$
 Thanks to $\omega_{n_j}\rightarrow\omega$ in $C([0,T_0];L^{p}(\R^2_+))$, it follows that  $\widetilde{\omega}=\omega$. Therefore, we have that
 $\omega\in L^\infty(0,T_0;W^{k,p}(\R^2_+))$
and
$\|\omega\|_{L^\infty(0,T_0;W^{k,p}(\R^2_+))}\leq M.$
Moreover, it is clear that $\omega\in C([0,T_0];L^{p}(\R^2_+))$ and $\omega(\cdot,0)=\omega_{0}\in W^{k,p}(\R^{2}_+).$
Therefore, $\omega\in B_{T_0},$
which shows that $B_{T_0}$ is a closed subset of $C([0,T_0];L^{p}(\R^2_+))$.\\[2mm]
{\it Step 3. The mapping $\mathcal{T}$ is  contractive on $B_{T_0}$ in the topology of $C([0,T_0];L^{p}(\R^2_+))$.}

Suppose that $\omega_i\in B_{T_0}$, $u_i$ is determined from $\omega_i$ by \eqref{3.3} and $\theta_i=\mathcal{T}\omega_i$ with $i=1,2.$
It follows from \eqref{4.2} that
\begin{equation}\label{4.5}
\left\{\ba
&\partial_t\theta_1+u_1\cdot\nabla\theta_1=0,\\
&\partial_t\theta_2+u_2\cdot\nabla\theta_2=0. \ea\ \right.
\end{equation}
Subtracting the equations in \eqref{4.5}, we get
\begin{eqnarray}\label{4.6}
\partial_t(\theta_1-\theta_2)+(u_1-u_2)\cdot\nabla\theta_1+u_2\cdot\nabla(\theta_1-\theta_2)=0.
\end{eqnarray}
Multiplying  by $|\theta_1-\theta_2|^{p-2}(\theta_1-\theta_2)$ on both sides of \eqref{4.6}, integrating the resulting equation on $\R^2_+$, making
integration by parts, applying  H\"{o}lder inequality and \eqref{3.4}-\eqref{3.5}, we  have
\begin{eqnarray}\label{4.7}
&&\frac{1}{p}\frac{d}{dt}\|\theta_1(t)-\theta_2(t)\|^{p}_{L^{p}(\R^2_+)}\nonumber\\
&=&-\int_{\mathbb{R}^2_+}((u_1-u_2)\cdot\nabla\theta_1)|\theta_1-\theta_2|^{p-2}(\theta_1-\theta_2)dx\nonumber\\
&\leq&\|u_1-u_2\|_{L^{\frac{2p}{2-p}}(\R^2_+)}\|\nabla\theta_1\|_{L^2(\R^2_+)}\|\theta_1-\theta_2\|^{p-1}_{L^{p}(\R^2_+)}\nonumber\\
&\leq&C\|\omega_1-\omega_2\|_{L^{p}(\R^2_+)}\|\theta_1\|_{W^{k,p}(\R^2_+)}\|\theta_1-\theta_2\|^{p-1}_{L^{p}(\R^2_+)}\nonumber\\
&\leq&CM\|\omega_1(t)-\omega_2(t)\|_{L^{p}(\R^2_+)}\|\theta_1(t)-\theta_2(t)\|^{p-1}_{L^{p}(\R^2_+)},
\end{eqnarray}
where we also use the inequality
\begin{eqnarray}\label{4.8}
\|f\|_{L^2(\R^2_+)}
&\leq&C\|f\|_{W^{2,p}(\R^2_+)},
\end{eqnarray}
for $1<p<2$, which can be deduced from the interpolation inequality in $L^s(\R^2_+)$ space and the Sobolev embedding $W^{2,s}(\R^{2}_+)\hookrightarrow L^\infty(\R^2_+)$ for any $s>1$.
It follows from \eqref{4.7} that
\begin{eqnarray*}
\frac{d}{dt}\|\theta_1(t)-\theta_2(t)\|_{L^{p}(\R^2_+)}\leq CM\|\omega_1(t)-\omega_2(t)\|_{L^{p}(\R^2_+)}
\end{eqnarray*}
Noticing that $\theta_1(\cdot,0)=\theta_2(\cdot,0)$, we obtain, for a.e. $t\in [0,T_0]$,
\begin{eqnarray*}
\|\theta_1(t)-\theta_2(t)\|_{L^{p}(\R^2_+)}
&\leq&\|\theta_1(0)-\theta_2(0)\|_{L^{p}(\R^2_+)}+CM\int_{0}^{t}\|\omega_1(\tau)-\omega_2(\tau)\|_{L^{p}(\R^2_+)}d\tau\\
&\leq& CM\int_{0}^{T_0}\|\omega_1(\tau)-\omega_2(\tau)\|_{L^{p}(\R^2_+)}d\tau\\
&\leq& CMT_0\|\omega_1-\omega_2\|_{C([0,T_0];L^{p}(\R^2_+))}.
\end{eqnarray*}
Since $\theta_1-\theta_2\in C([0,T_0];L^{p}(\R^2_+))$, we then have
\begin{eqnarray*}
\|\theta_1-\theta_2\|_{C([0,T_0];L^{p}(\R^2_+))}\leq CMT_0\|\omega_1-\omega_2\|_{C([0,T_0];L^{p}(\R^2_+))}.
\end{eqnarray*}
Choosing $T_0=\min\{\frac{\log2}{CM},\frac{1}{2CM}\}$,   we obtain
\begin{eqnarray*}
\|\mathcal{T}\omega_1-\mathcal{T}\omega_2\|_{C([0,T_0];L^{p}(\R^2_+))}\leq \frac{1}{2}\|\omega_1-\omega_2\|_{C([0,T_0];L^{p}(\R^2_+))},
\end{eqnarray*}
which implies that $\mathcal{T}$ is contractive on $ B_{T_0}$ in the topology of $C([0,T_0];L^{p}(\R^2_+))$.
Finally, combining step 1, step 2 and step 3, by virtue of Lemma \ref{contration}, we finish the proof of Lemma \ref{existence of fixed-point}.

Now we are ready to prove Theorem \ref{the-1}.

\textbf{Proof of Theorem \ref{the-1}.}
By virtue of Lemma \ref{existence of fixed-point}, there exists exactly one $\omega\in B_{T_0}$ such that $\theta=\mathcal{T}\omega=\omega,$ which implies that
\begin{equation*}
\left\{\ba
&\partial_t\omega+u\cdot\nabla\omega=0, ~(x,t)\in \R^{2}_+\times\R_+\\&u =\nabla^\bot\Delta^{-1}_D\omega\\
&\omega(x,0)=\omega_{0}(x), \ea\ \right.
\end{equation*}
Moreover, in view of Lemma \ref{Linear-Transport}, the solution $\omega$ belongs to $L^\infty(0,T_0;W^{k,p}(\R^2_+))\cap{\rm Lip}([0,T_0];W^{k-1,p}(\R^2_+)).$

The uniqueness is sketched as follows, which is similar to {\it Step 3} in the proof of Lemma \ref{existence of fixed-point}.
Suppose that $\omega_i\in L^\infty(0,T_0;W^{k,p}(\R^2_+))\cap{\rm Lip}([0,T_0];W^{k-1,p}(\R^2_+))$ with $k\geq3$ and $1<p<2$, $i=1,2,$ are two solutions to \eqref{1.4} with the same initial data $\omega_0\in W^{k,p}(\R^2_+)$. Then
\begin{equation}\label{4.9}
\left\{\ba
&\partial_t\omega_1+u_1\cdot\nabla\omega_1=0,\\
&\partial_t\omega_2+u_2\cdot\nabla\omega_2=0. \ea\ \right.
\end{equation}
Subtracting the equations in \eqref{4.9} yields
\begin{eqnarray*}
\partial_t(\omega_1-\omega_2)+(u_1-u_2)\cdot\nabla\omega_1+u_2\cdot\nabla(\omega_1-\omega_2)=0.
\end{eqnarray*}
Similar to \eqref{4.7}, we can obtain
\begin{eqnarray*}
\frac{d}{dt}\|\omega_1(t)-\omega_2(t)\|_{L^{p}(\R^2_+)}
\leq C\|\omega_1(t)\|_{W^{k,p}(\R^2_+)}\|\omega_1(t)-\omega_2(t)\|_{L^{p}(\R^2_+)}.
\end{eqnarray*}
Applying Gronwall's inequality leads to
$$\|\omega_1(t)-\omega_2(t)\|_{L^{p}(\R^2_+)}
\leq e^{C\int_{0}^{t}\|\omega_1(\tau)\|_{W^{k,p}(\R^2_+)}d\tau}\|\omega_1(0)-\omega_2(0)\|_{L^{p}(\R^2_+)},$$
which implies that $\omega_1=\omega_2$.

Finally, we verify that $\omega\in C([0,T_0];W^{k,p}(\R^2_+)),$ which is simialr to the case of whole plane (see, e.g.,  \cite{[MB]}). In fact, it is direct to deduce that $\omega\in C_w([0,T_0];W^{k,p}(\R^2_+))$, which denotes continuity on the interval $[0,T_0]$ with values in the weak topology of $W^{k,p}(\R^2_+)$.
To prove the strong continuity in $\omega\in C([0,T_0];W^{k,p}(\R^2_+))$, it suffices to show that the norm function $\|\omega(t)\|_{W^{k,p}(\R^2_+)}$
is continuous in time. We first  prove the continuity of $\|\omega(t)\|_{W^{k,p}(\R^2_+)}$ at the initial time.

On one hand, by the weak continuity of the resulting solution, we have $$\|\omega_0\|_{W^{k,p}(\R^2_+)}\leq \displaystyle\liminf_{t\rightarrow 0^+}\|\omega(t)\|_{W^{k,p}(\R^2_+)}.$$
On the other hand, since $\omega\in B_{T_0}$,   it follows  from \eqref{4.3} and \eqref{4.4} that, for a.e. $t\in[0,T_0]$,
\begin{eqnarray*}
\|\mathcal{T}\omega(t)\|_{W^{k,p}(\R^2_+)}&\leq&\|\omega_0\|_{W^{k,p}(\R^2_+)}e^{C\int_{0}^{t}\|\omega(\tau)\|_{W^{k,p}(\R^2_+)}d\tau}\\
&\leq&\|\omega_0\|_{W^{k,p}(\R^2_+)}e^{CMt}.
\end{eqnarray*}
Since $\omega$ is a fixed point of $\mathcal{T}$, we obtain, for a.e. $t\in[0,T_0]$,
\begin{eqnarray*}
\|\omega(t)\|_{W^{k,p}(\R^2_+)}
&\leq&\|\omega_0\|_{W^{k,p}(\R^2_+)}e^{CMt}.
\end{eqnarray*}
Consequently,
\begin{eqnarray*}
\displaystyle\limsup_{t\rightarrow 0^+}\|\omega(t)\|_{W^{k,p}(\R^2_+)}\leq \|\omega_0\|_{W^{k,p}(\R^2_+)}.
\end{eqnarray*}
Hence it holds that
\begin{eqnarray*}
\displaystyle\lim_{t\rightarrow 0^+}\|\omega(t)\|_{W^{k,p}(\R^2_+)}=\|\omega_0\|_{W^{k,p}(\R^2_+)},
\end{eqnarray*}
and the solution $\omega$ is strongly right continuous at $t=0.$
To prove the same strong continuity of the solution $\omega$ at any time $t_0\in(0,T_0]$, let $\widehat{\omega}$ be the local solution to \eqref{1.4} for $t\geq t_0$ with the initial value $\omega_0(t_0)$. By the result just proved, $\|\widehat{\omega}(t)\|_{W^{k,p}(\R^2_+)}$ is right continuous at $t=t_0.$ But $\omega$ coincides with $\widehat{\omega}$ for $t\geq t_0 $ by the uniqueness. Hence $\|\omega(t)\|_{W^{k,p}(\R^2_+)}$ is right continuous at $t=t_0$. Due to the arbitrariness of $t_0\in(0,T_0]$, it follows that $\|\omega(t)\|_{W^{k,p}(\R^2_+)}$ is right continuous on  $(0,T_0]$.
Furthermore, since Euler equation is reversible in time $t$,  $\|\omega(t)\|_{W^{k,p}(\R^2_+)}$ must be also left continuous at $t=t_0$. Consequently, $\|\omega(t)\|_{W^{k,p}(\R^2_+)}$
is continuous on $[0,T_0]$ and $\omega\in C([0,T_0],W^{k,p}(\R^2_+))$.
The proof of Theorem \ref{the-1} is complete.
\section{Proof of Theorem \ref{the-2}}
In this section, we prove Theorem \ref{the-2} which is on global well-posedness to \eqref{1.4}. Before that, we   prove some useful facts for preparations. The first lemma is the maximum principle of the solutions to \eqref{1.4}, which states that the $L^s(\R^2_+)$ norm of the vorticity is conserved for all times.
\begin{lemma}\label{conserved quantity}
Let $1\leq s\leq\infty$. If $\omega$ is a smooth solution to \eqref{1.4}, then $\|\omega(t)\|_{L^s(\R^2_+)}=\|\omega_0\|_{L^s(\R^2_+)}$ for any $t>0$.
\end{lemma}
\textbf{Proof}.
Denote by $\Phi_t(x)$ the flow map corresponding to the two-dimentional Euler evolution in the half plane:
\begin{equation*}
\left\{\ba
&\frac{d}{dt}\Phi_t(x)=u(\Phi_t(x),x),\\\
&\Phi_0(x)=x.\ea\ \right.
\end{equation*}
Then it follows from \eqref{3.4} and sufficient regularity of $u$ that (see, e.g.,  \cite{[MB]})
\begin{eqnarray}\label{5.1}
\omega(x,t)=\omega_0(\Phi_{-t}(x))~{\rm and}~ {\rm det}\Big(\nabla_x\Phi_t(x)\Big)=1,
\end{eqnarray}
where $\Phi_{-t}$ is the inverse of $\Phi_{t}$. For $s\in[1,\infty)$, in view of \eqref{2.1}, we have
\begin{eqnarray*}
\|\omega(t)\|^s_{L^s(\R^2_+)}
&=&s\int^\infty_0\tau^{s-1}\mu_{\omega_t}(\tau)d\tau\\
&=&s\int^\infty_0\tau^{s-1}\mu_{\omega_0}(\tau)d\tau\\
&=&\|\omega_0\|^s_{L^s(\R^2_+)},
\end{eqnarray*}
where we have used the fact that the distribution function of the vorticity $\omega(\cdot,t):=\omega_t$
keeps same for all times:
\begin{eqnarray*}
\mu_{\omega_t}(\tau)=\mu_{\omega_0}(\tau)~ for~ all~ \tau>0.
\end{eqnarray*}
Indeed, by virtue of \eqref{2.2}, \eqref{5.1} and a change of variables, we obtain
\begin{eqnarray*}
\mu_{\omega_t}(\tau)
&=&|\{x\in\R^2_+:|\omega(x,t)|>\tau\}|\\
&=&|\{x\in\R^2_+:|\omega_0(\Phi_{-t}(x))|>\tau\}|\\
&=&\int_{\{x\in\R^2_+:|\omega_0(\Phi_{-t}(x))|>\tau\}}1~dx\\
&=&\int_{\{y\in\R^2_+:|\omega_0(y)|>\tau\}}{\rm det} \Big(\nabla_y\Phi_t(y)\Big)dy\\
&=&\int_{\{y\in\R^2_+:|\omega_0(y)|>\tau\}}1~dy\\
&=&|\{y\in\R^2_+:|\omega_0(y)|>\tau\}|\\
&=&\mu_{\omega_0}(\tau).
\end{eqnarray*}
For $s=\infty$, $\|\omega(t)\|_{L^\infty(\R^2_+)}=\|\omega_0\|_{L^\infty(\R^2_+)}$ is a consequence of \eqref{5.1}. The proof of the lemma is finished.

The second lemma is concerned with a Kato-type estimate of the velocity in the half plane case. The similar bound was first obtained by Kato \cite{[Kato2]} in the whole plane case. We sketch its proof for completeness in the spirit of \cite{[KS]}, where the bounded domain case was treated. The similar estimates in three-dimensional whole space and smooth bounded domains of $\R^3$ can be found in \cite{[BKM]} and \cite{[F]}, respectively.
\begin{lemma}\label{Kato inequality}
For any $0<\gamma<1$ and $1\leq q<2$, it holds
\begin{eqnarray*}
\| u(t)\|_{W^{1,\infty}(\R^2_+)}\leq C_\gamma\|\omega_0\|_{L^\infty(\R^2_+)}\Big(1+\log\Big(1+\frac{[\omega(t)]_{C^\gamma(\R^2_+)}}{\|\omega_0\|_{L^\infty(\R^2_+)}}\Big)\Big)
+C_q\|\omega_0\|_{L^q(\R^2_+)}.
\end{eqnarray*}
Consequently,
\begin{eqnarray*}
\|u(t)\|_{W^{1,\infty}(\R^2_+)}\leq C_q(\|\omega_0\|_{L^\infty(\R^2_+)}+\|\omega_0\|_{L^q(\R^2_+)})
\Big(1+\log\Big(3+\frac{\|\nabla\omega(t)\|_{L^\infty(\R^2_+)}}{\|\omega_0\|_{L^\infty(\R^2_+)}}\Big)\Big)
,
\end{eqnarray*}
for $1\leq q<2$.
\end{lemma}
\textbf{Proof}. We first estimate $\|u(t)\|_{L^\infty(\R^2_+)}$.
From \eqref{3.3}, the velociry is expressed as
\begin{eqnarray*}
u(x,t)
&=&\frac{1}{2\pi}\int_{\mathbb R^2_+}\Big(\frac{(x-y)^\bot}{|x-y|^2}-\frac{(x-\bar{y})^\bot}{|x-\bar{y}|^2}\Big)\omega(y,t)dy\\
&=&\frac{1}{2\pi}\int_{\mathbb R^2}\frac{(x-y)^\bot}{|x-y|^2}\bar{\omega}(y,t)dy,
\end{eqnarray*}
where $\bar{\omega}:\R^2\rightarrow \R $ is the odd extension of $\omega$ to the whole plane.
Therefore, for any $x\in \R^2_+$ and $1\leq q<2$,  by virtue of H\"{o}lder inequality and Lemma \ref{conserved quantity}, it yields
\begin{eqnarray*}
|u(x,t)|
&\leq&\frac{1}{2\pi}\int_{\mathbb R^2}\frac{1}{|x-y|}|\bar{\omega}(y,t)|dy\\
&=&\frac{1}{2\pi}\int_{|x-y|\leq1}\frac{1}{|x-y|}|\bar{\omega}(y,t)|dy
+\frac{1}{2\pi}\int_{|x-y|>1}\frac{1}{|x-y|}|\bar{\omega}(y,t)|dy\\
&\leq&\frac{1}{2\pi}\|\bar{\omega}(t)\|_{L^\infty(\R^2)}\int_{|x-y|\leq1}\frac{dy}{|x-y|}
+\frac{1}{2\pi}\|\bar{\omega}(t)\|_{L^q(\R^2)}\||z|^{-1}\|_{L^{q\prime}((B_1(x))^c,dz)}\\
&\leq&\|\omega(t)\|_{L^\infty(\R^2_+)}
+C_q\|\omega(t)\|_{L^q(\R^2_+)}\\
&=&\|\omega_0\|_{L^\infty(\R^2_+)}+C_q\|\omega_0\|_{L^q(\R^2_+)},
\end{eqnarray*}
where $q'=\frac{q}{q-1}$ for $1<q<2$ and we define $q'=\infty$ when $q=1$.
It shows that
\begin{eqnarray*}
\|u(t)\|_{L^\infty(\R^2_+)}\leq\|\omega_0\|_{L^\infty(\R^2_+)}+C_q\|\omega_0\|_{L^q(\R^2_+)}.
\end{eqnarray*}
Hence it remains to estimate $\|\nabla u\|_{L^\infty(\R^2_+)}$.
By Lemma \ref{Formula for gradient}, the gradient of the velocity is expressed as
\begin{eqnarray}\label{5.2}
\nabla u(x,t)=\frac{1}{2\pi}P.V.\int_{\R^2_+}M(x,y)\omega(y,t)dy+\frac{\omega(x,t)}{2}
\begin{pmatrix}
0&1\\-1&0
\end{pmatrix},
\end{eqnarray}
where the kernel matrix $M(x,y)$ is the sum of singular  part $M_s(x,y)$ and regular part $M_r(x,y)$:
\begin{eqnarray*}
M_s(x,y)=
\begin{pmatrix}
\frac{-2(x_1-y_1)(x_2-y_2)}{|x-y|^4}
&\frac{(x_1-y_1)^2-(x_2-y_2)^2}{|x-y|^4}\\
\frac{(x_1-y_1)^2-(x_2-y_2)^2}{|x-y|^4}
&\frac{2(x_1-y_1)(x_2-y_2)}{|x-y|^4}
\end{pmatrix},
\end{eqnarray*}
\begin{eqnarray*}
M_r(x,y)=
\begin{pmatrix}
\frac{2(x_1-y_1)(x_2+y_2)}{|x-\bar{y}|^4}
&-\frac{(x_1-y_1)^2-(x_2+y_2)^2}{|x-\bar{y}|^4}\\
-\frac{(x_1-y_1)^2-(x_2+y_2)^2}{|x-\bar{y}|^4}
&-\frac{2(x_1-y_1)(x_2+y_2)}{|x-\bar{y}|^4}
\end{pmatrix}.
\end{eqnarray*}
The second term of \eqref{5.2} is controlled by $C\|\omega(t)\|_{L^\infty(\R^2_+)}\leq C\|\omega_0\|_{L^\infty(\R^2_+)}$ by Lemma \ref{conserved quantity}.
Therefore it suffices to deal with the principle integral in \eqref{5.2}. We proceed to set $\delta(t):=\Big(\frac{\|\omega(t)\|_{L^\infty(\R^2_+)}}{\|\omega(t)\|_{C^\gamma(\R^2_+)}}\Big)^{\frac{1}{\gamma}}$. Then it is clear that $0<\delta(t)\leq1$.

For an interior point $x=(x_1,x_2)\in\R^2_+$ with $x_2={\rm dist}(x,\partial\R^2_+)> 2\delta(t)$, we have $B_{2\delta(t)}(x)\subseteq\R^2_+$. Then it follows that
\begin{eqnarray}\label{5.3}
|P.V.\int_{\R^2_+}M(x,y)\omega(y,t)dy|\leq I_1+I_2+I_3,
\end{eqnarray}
where
\begin{eqnarray}\label{5.4}
I_1
&=&|P.V.\int_{\R^2_+\cap B_{\delta(t)}(x)}M(x,y)\omega(y,t)dy|\nonumber\\
&\leq&|P.V.\int_{\R^2_+\cap B_{\delta(t)}(x)}M_s(x,y)\omega(y,t)dy|
+|\int_{\R^2_+\cap B_{\delta(t)}(x)}M_r(x,y)\omega(y,t)dy|\nonumber\\
&=:& I_{11}+I_{12},
\end{eqnarray}
\begin{eqnarray*}
I_2=|\int_{\R^2_+\cap (B_1(x)\setminus B_{\delta(t)}(x))}M(x,y)\omega(y,t)dy|,
\end{eqnarray*}
and
\begin{eqnarray*}
I_3=|\int_{\R^2_+\cap(B_1(x))^c}M(x,y)\omega(y,t)dy|.
\end{eqnarray*}
For $I_{11}$, utilizing $P.V.\int_{|x-y|<\delta(t)}M_s(x,y)dy=0$, $|M_s(x,y)|\leq C|x-y|^{-2}$, Lemma \ref{5.1} and the definition of $\delta(t)$, we have
\begin{eqnarray}\label{5.5}
I_{11}
&=&|P.V.\int_{|x-y|<\delta(t)}M_s(x,y)\Big(\omega(y,t)-\omega(x,t)\Big)dy|\nonumber\\
&\leq&P.V.\int_{|x-y|<\delta(t)}|M_s(x,y)||\omega(y,t)-\omega(x,t)|dy\nonumber\\
&\leq&C[\omega(t)]_{C^\gamma(\R^2_+)}\int_{|x-y|<\delta(t)}\frac{1}{|x-y|^2}|x-y|^\gamma dy\nonumber\\
&\leq&C\|\omega(t)\|_{C^\gamma(\R^2_+)}\int_{|z|<\delta(t)}\frac{dz}{|z|^{2-\gamma}}\nonumber\\
&=&\frac{C}{\gamma}(\delta(t))^\gamma\|\omega(t)\|_{C^\gamma\R^2_+)}\nonumber\\
&\leq&\frac{C}{\gamma}\|\omega(t)\|_{L^\infty(\R^2_+)}\nonumber\\
&=&\frac{C}{\gamma}\|\omega_0\|_{L^\infty(\R^2_+)}.
\end{eqnarray}
For $I_{12}$, using $|M_r(x,y)|\leq C|x-\bar{y}|^{-2}$, Lemma \ref{5.1} and $x_2>2\delta(t)$, we can estimate $I_{12}$ as follows.
\begin{eqnarray}\label{5.6}
I_{12}
&\leq&\int_{\R^2_+\cap B_{\delta(t)}(x)}|M_r(x,y)||\omega(y,t)|dy\nonumber\\
&\leq&C\int_{\R^2_+\cap B_{\delta(t)}(x)}|x-\bar{y}|^{-2}|\omega(y,t)|dy\nonumber\\
&\leq&C\|\omega(t)\|_{L^\infty(\R^2_+)}\Big(\frac{\delta(t)}{x_2}\Big)^2\nonumber\\
&\leq&C\|\omega_0\|_{L^\infty(\R^2_+)},
\end{eqnarray}
where we have used $|x-\bar{y}|>x_2$ for any $x, y\in\R^2_+$.

For $I_2$, utilizing $|M(x,y)|\leq C|x-y|^{-2}$, H\"{o}lder inequality, the definition of $\delta(t)$ and Lemma \ref{conserved quantity}, we obtain
\begin{eqnarray}\label{5.7}
I_2
&\leq&\int_{\R^2_+\cap (B_1(x)\setminus B_{\delta(t)}(x))}|M(x,y)||\omega(y,t)|dy\nonumber\\
&\leq&C\|\omega(t)\|_{L^\infty(\R^2_+)}\int_{\delta(t)\leq|x-y|<1}\frac{dy}{|x-y|^2}\nonumber\\
&\leq&C\|\omega(t)\|_{L^\infty(\R^2_+)}\int_{\delta(t)\leq|z|<1}\frac{dz}{|z|^2}\nonumber\\
&=&C\|\omega(t)\|_{L^\infty(\R^2_+)}\log\frac{1}{\delta(t)}\nonumber\\
&=&\frac{C}{\gamma}\|\omega(t)\|_{L^\infty(\R^2_+)}\log\Big(1+\frac{[\omega(t)]_{C^\gamma(\R^2_+)}}{\|\omega(t)\|_{L^\infty(\R^2_+)}}\Big)\nonumber\\
&\leq&\frac{C}{\gamma}\|\omega_0\|_{L^\infty(\R^2_+)}\log\Big(1+\frac{[\omega(t)]_{C^\gamma(\R^2_+)}}{\|\omega_0\|_{L^\infty(\R^2_+)}}\Big),
\end{eqnarray}
where in the last inequality we have used the  fact that the function $x\mapsto x{\rm log}(1+\frac{\beta}{x})$ is non-decreasing on $(0,\infty)$ for all $\beta\geq0$.

For $I_3$, using $|M(x,y)|\leq C|x-y|^{-2}$, H\"{o}lder inequality and Lemma \ref{conserved quantity}, we obtain, for $1\leq q<2$,
\begin{eqnarray}\label{5.8}
I_3
&\leq&\int_{\R^2_+\cap(B_1(x))^c}|M(x,y)||\omega(y,t)|dy\nonumber\\
&\leq&C\int_{\R^2_+\cap(B_1(x))^c}\frac{|\omega(y,t)|}{|x-y|^2}dy\nonumber\\
&\leq&C\|\omega(t)\|_{L^q(\R^2_+)}\||z|^{-2}\|_{L^{q'}((B_1(x))^c,dz)}\nonumber\\
&\leq&C_q\|\omega(t)\|_{L^q(\R^2_+)}\nonumber\\
&=&C_q\|\omega_0\|_{L^q(\R^2_+)}.
\end{eqnarray}
Combining \eqref{5.3}-\eqref{5.8} yields the desired estimate of $\nabla u$ at interior points.

 For a point nearby the boundary, $x=(x_1,x_2)\in\R^2_+$ with $x_2={\rm dist}(x,\partial\R^2_+)\leq 2\delta(t)$, we define $z=(x_1,x_2+2\delta(t))$. Then it is clear that ${\rm dist}(z,\partial\R^2_+)=x_2+2\delta(t)> 2\delta(t)$ and $|x-z|=2\delta(t)$.
It follows that
\begin{eqnarray*}
|\nabla u(x,t)|
&\leq&|\nabla u(x,t)-\nabla u(z,t)|+|\nabla u(z,t)|\\
&\leq&[\nabla u(t)]_{C^\gamma(\R^2_+)}|x-z|^\gamma+|\nabla u(z,t)|\\
&=&2^\gamma(\delta(t))^\gamma[\nabla u(t)]_{C^\gamma(\R^2_+)}+|\nabla u(z,t)|.
\end{eqnarray*}
Utilizing Proposition \ref{Schauder-type-estimate}, the definition of $\delta(t)$, Lemma \ref{conserved quantity} and the  estimate of $\nabla u$ at interior points $z$ above, we  obtain, for $0<\gamma<1$ and $1\leq q<2$,
\begin{eqnarray*}
|\nabla u(x,t)|
&\leq&C_\gamma(\delta(t))^\gamma[\omega(t)]_{C^\gamma(\R^2_+)}+|\nabla u(z,t)|\\
&\leq&C_\gamma\|\omega_0\|_{L^\infty(\R^2_+)}\Big(1+\log\Big(1+\frac{[\omega(t)]_{C^\gamma(\R^2_+)}}{\|\omega_0\|_{L^\infty(\R^2_+)}}\Big)\Big)
+C_q\|\omega_0\|_{L^q(\R^2_+)},
\end{eqnarray*}
which is the desired estimate of $\nabla u$ at points near the boundary. Furthermore, applying the fact $[\omega(t)]_{C^\gamma(\R^2_+)}\leq2\|\omega_0\|_{L^\infty(\R^2_+)}+\|\nabla\omega(t)\|_{L^\infty(\R^2_+)}$ finally implies that
\begin{eqnarray*}
\|u(t)\|_{W^{1,\infty}(\R^2_+)}\leq C_q(\|\omega_0\|_{L^\infty(\R^2_+)}+\|\omega_0\|_{L^q(\R^2_+)})
\Big(1+\log\Big(3+\frac{\|\nabla\omega(t)\|_{L^\infty(\R^2_+)}}{\|\omega_0\|_{L^\infty(\R^2_+)}}\Big)\Big)
,
\end{eqnarray*}
for $1\leq q<2$.

The proof of the lemma  is complete.

We are now in a position to derive the double exponential growth in time of the gradient vorticity of the two-dimensional incompressible Euler equations in the half plane, which is a key ingredient in the proof of global regularity of the solutions to \eqref{1.4}.
\begin{proposition}\label{Double-Exponetial-Upper-Bound}
Let $\omega_0$ be smooth initial data for the 2-D Euler equation \eqref{1.4} in the half plane. Then the solution $\omega$ to \eqref{1.4} satisfies
\begin{eqnarray*}
1+\log\Big(3+\frac{\|\nabla\omega(t)\|_{L^\infty(\R^2_+)}}{\|\omega_0\|_{L^\infty(\R^2_+)}}\Big)
\leq \Big(1+\log\Big(3+\frac{\|\nabla\omega_0\|_{L^\infty(\R^2_+)}}{\|\omega_0\|_{L^\infty(\R^2_+)}}\Big)\Big)
e^{C(\|\omega_0\|_{L^\infty(\R^2_+)}+\|\omega_0\|_{L^{q}(\R^2_+)})t},
\end{eqnarray*}
where the constant $C$ depends on $q\in[1,2)$.
\end{proposition}
\textbf{Proof}.
 Applying the gradient operator $\nabla$ on the first equation in \eqref{1.4}, and taking its dot product by $\nabla \omega$, we obtain
\begin{eqnarray*}
\partial_t\Big(\frac{1}{2}|\nabla\omega|^2\Big)+u\cdot\nabla\Big(\frac{1}{2}|\nabla\omega|^2\Big)+\sum_{1\leq i,j\leq2} \partial_iu^j\partial_i\omega\partial_j\omega=0.
\end{eqnarray*}
It follows that
\begin{eqnarray}\label{5.9}
\partial_t|\nabla\omega|+u\cdot\nabla|\nabla\omega|=-\beta|\nabla\omega|,
\end{eqnarray}
where
\begin{eqnarray*}
\beta(x,t):=\sum_{1\leq i,j\leq2} \partial_iu^j\frac{\partial_i\omega}{|\nabla\omega|}\frac{\partial_j\omega}{|\nabla\omega|}.
\end{eqnarray*}
Then
\begin{eqnarray*}
\frac{d}{dt}|\nabla\omega|(\Phi_t(x),t)=-\beta(\Phi_t(x),t)|\nabla\omega|(\Phi_t(x),t),
\end{eqnarray*}
where $\Phi$ is the flow map  generated by the velocity field $u$. Hence,
\begin{eqnarray*}
|\nabla\omega|(\Phi_t(x),t)
&=&|\nabla\omega_0(x)|e^{-\int_{0}^{t}\beta(\Phi_\tau(x),\tau)d\tau}\\
&\leq&
|\nabla\omega_0(x)|e^{\int_{0}^{t}|\nabla u|(\Phi_\tau(x),\tau)d\tau},
\end{eqnarray*}
where we have used $|\beta(x,t)|\leq|\nabla u(x,t)|$ in the last inequality.
It deduces that
\begin{eqnarray*}
\|\nabla\omega(t)\|_{L^\infty(\R^2_+)}
&\leq&
\|\nabla\omega_0\|_{L^\infty(\R^2_+)}e^{\int_{0}^{t}\|\nabla u(\tau)\|_{L^\infty(\R^2_+)}d\tau}.
\end{eqnarray*}
It follows from Lemma \ref{Kato inequality} that
\begin{eqnarray*}
\log\Big(\frac{\|\nabla\omega(t)\|_{L^\infty(\R^2_+)}}{\|\omega_0\|_{L^\infty(\R^2_+)}}\Big)
\leq \log\Big(\frac{\|\nabla\omega_0\|_{L^\infty(\R^2_+)}}{\|\omega_0\|_{L^\infty(\R^2_+)}}\Big) +CA\int_{0}^{t}\Big(1+\log\Big(3+\frac{\|\nabla\omega(\tau)\|_{L^\infty(\R^2_+)}}{\|\omega_0\|_{L^\infty(\R^2_+)}}\Big)\Big)d\tau,
\end{eqnarray*}
with $A:=\|\omega_0\|_{L^\infty(\R^2_+)}+\|\omega_0\|_{L^q(\R^2_+)}$.
Consider
\begin{equation}\label{5.10}
\left\{\ba
&\frac{d}{dt}\log z(t)=CA(1+\log(3+z(t))),\\
&z(0)=\frac{\|\nabla\omega_0\|_{L^\infty(\R^2)}}{\|\omega_0\|_{L^\infty(\R^2)}}.\ea\ \right.
\end{equation}
By Lemma \ref{Integral-Inequality}, we have$\frac{\|\nabla\omega(t)\|_{L^\infty(\R^2_+)}}{\|\omega_0\|_{L^\infty(\R^2_+)}}\leq z(t)$. To bound $z(t)$, we solve \eqref{5.10} directly to obtain
\begin{eqnarray*}
\int^{z(t)}_{z(0)}\frac{dy}{y(1+\log(3+y))}=CAt.
\end{eqnarray*}
Therefore,
\begin{eqnarray*}
\int^{z(t)}_{z(0)}\frac{dy}{(3+y)(1+\log(3+y))}\leq CAt,
\end{eqnarray*}
and
\begin{eqnarray*}
1+\log\Big(3+z(t)\Big)
\leq \Big(1+\log\Big(3+z(0)\Big)\Big)
\exp(CAt),
\end{eqnarray*}
which concludes the proof of Lemma \ref{Double-Exponetial-Upper-Bound}.

Combining Lemma \ref{Kato inequality} and Proposition \ref{Double-Exponetial-Upper-Bound} immediately yields the estimate on $\|u\|_{L^\infty(0,T;W^{1,\infty}(\R^2_+))}$ for any $T>0$.
\begin{coro}\label{single-exponential-velocity-gradient}
For any $1\leq q<2$, it holds that, for any $t>0$,
\begin{eqnarray*}
\|u(t)\|_{W^{1,\infty}(\R^2_+)}\leq C_qABe^{C_qAt},
\end{eqnarray*}
with
\begin{eqnarray*}
A=\|\omega_0\|_{L^\infty(\R^2_+)}+\|\omega_0\|_{L^q(\R^2_+)},
\end{eqnarray*}
and
\begin{eqnarray*}
B=1+\log\Big(3+\frac{\|\nabla\omega_0\|_{L^\infty(\R^2_+)}}{\|\omega_0\|_{L^\infty(\R^2_+)}}\Big).
\end{eqnarray*}
\end{coro}
Now we are ready to prove Theorem \ref{the-2}.

\textbf{Proof of Theorem \ref{the-2}.}
For any fixed $T>0$. By Theorem \ref{the-1}, \eqref{1.4} has a local unique solution $\omega\in C([0,T^\ast);W^{k,p}(\R^{2}_+))$ for the initial data $\omega_0\in W^{k,p}(\R^{2}_+)$ with $k\geq3$ and $1<p<2$ , where $T^\ast$ is the first blow-up time.
To prove the global regularity of $\omega$, it suffices to show that $T^\ast>T$. By the classical bootstrap and continuity argument, we only need to derive the uniform estimate of $\displaystyle\sup_{0\leq t<T^\ast}\|\omega(t)\|_{W^{k,p}(\R^2_+)}$.
  For this purpose, by a standard process of approximating the initial data and local well-posedness Theorem \ref{the-1}, we can assume $\omega\in C([0,T^\ast);W^{k+1,p}(\R^2_+))$. Then,
for arbitrary multi-indice $\alpha\in\mathbb{N}^{2}$ with $0<|\alpha|\leq k$, applying $\partial^{\alpha}$ on the first equation of \eqref{1.4} gives
\begin{eqnarray}\label{5.11}
\partial_t\partial^{\alpha}\omega+u\cdot\nabla \partial^{\alpha}\omega+\partial^{\alpha}(u\cdot\nabla\omega)-u\cdot\nabla \partial^{\alpha}\omega=0.
\end{eqnarray}
Multiplying \eqref{5.11} by $|\partial^{\alpha}\omega|^{p-2}\partial^{\alpha}\omega$, integrating the resulting equation on $\R^2_+$, making
integration by parts, and utilizing \eqref{3.4}, H\"{o}lder inequality, Corollary \ref{Calculus-half-space} and Proposition \ref{Sobolev-type estimate}, we can finally arrive at
\begin{eqnarray*}
&&\frac{1}{p}\frac{d}{dt}\|\partial^{\alpha}\omega(t)\|^{p}_{L^{p}(\R^2_+)}\\
&=&-\int_{\mathbb{R}^2_+}(\partial^{\alpha}(u\cdot\nabla\omega)-u\cdot\nabla \partial^{\alpha}\omega)|\partial^{\alpha}\omega|^{p-2}(\partial^{\alpha}\omega)dx\\
&\leq&\|\partial^{\alpha}(u\cdot\nabla\omega)-u\cdot\nabla \partial^{\alpha}\omega\|_{L^{p}(\R^2_+)}\|\partial^{\alpha}\omega\|^{p-1}_{L^{p}(\R^2_+)}\\
&\leq&C\Big(\|u\|_{W^{|\alpha|,\frac{2p}{2-p}}(\R^2_+)}\|\nabla\omega\|_{L^2(\R^2_+)}
+\|u\|_{W^{1,\infty}(\R^2_+)}\|\nabla\omega\|_{W^{|\alpha|-1,p}(\R^2_+)}\Big)\|\partial^{\alpha}\omega\|^{p-1}_{L^{p}(\R^2)}\\
&\leq&C\Big(\|\nabla \omega\|_{L^2(\R^2_+)}+\|u\|_{W^{1,\infty}(\R^2_+)}\Big)\|\omega(t)\|_{W^{k,p}(\R^2_+)}\|\partial^{\alpha}\omega\|^{p-1}_{L^{p}(\R^2_+)}.
\end{eqnarray*}
It follows that, for any $0<|\alpha|\leq k$,
\begin{eqnarray}\label{5.12}
\frac{d}{dt}\|\partial^{\alpha}\omega(t)\|_{L^{p}(\R^2_+)}
\leq C\Big(\|\nabla \omega\|_{L^2(\R^2_+)}+\|u\|_{W^{1,\infty}(\R^2_+)}\Big)\|\omega(t)\|_{W^{k,p}(\R^2_+)}.
\end{eqnarray}
Moreover, it is straightforward to obtain
\begin{eqnarray}\label{5.13}
\frac{d}{dt}\|\omega(t)\|_{L^{p}(\R^2_+)}=0.
\end{eqnarray}
Combining \eqref{5.12} with \eqref{5.13} implies that
\begin{eqnarray}\label{5.14}
\frac{d}{dt}\|\omega(t)\|_{W^{k,p}(\R^2_+)}
\leq C\Big(\|\nabla \omega\|_{L^2(\R^2_+)}+\|u\|_{W^{1,\infty}(\R^2_+)}\Big)\|\omega(t)\|_{W^{k,p}(\R^2_+)}.
\end{eqnarray}
Hence, applying the Gronwall's inequality leads to
\begin{eqnarray}\label{5.15}
\|\omega(t)\|_{W^{k,p}(\R^2_+)}
\leq \|\omega_0\|_{W^{k,p}(\R^2_+)}e^{C{\int_{0}^{t}\|\nabla \omega(\tau)\|_{L^2(\R^2_+)}+\| u(\tau)\|_{W^{1,\infty}(\R^2_+)}d\tau}},
\end{eqnarray}
for all $t\in[0,T^\ast)$.

We proceed to give the  estimate of $\|\nabla\omega\|_{L^2(\R^2_+)}$.
Taking the inner product of \eqref{5.9} with $|\nabla \omega|$  in $L^2(\R^2_+)$, making integration by parts, utilizing \eqref{3.4} and H\"{o}lder inequality, we can obtain
\begin{eqnarray*}
\frac{1}{2}\frac{d}{dt}\|\nabla\omega(t)\|^2_{L^2(\R^2_+)}
&=&-\int_{\R^2_+}\beta|\nabla\omega|^2dx\\
&\leq&\|\beta(t)\|_{L^\infty(\R^2_+)}\|\nabla\omega(t)\|^2_{L^2(\R^2_+)}\\
&\leq&\|\nabla u(t)\|_{L^\infty(\R^2_+)}\|\nabla\omega(t)\|^2_{L^2(\R^2_+)},
\end{eqnarray*}
which deduces that
\begin{eqnarray*}
\frac{d}{dt}\|\nabla\omega(t)\|_{L^2(\R^2_+)}\leq\|\nabla u(t)\|_{L^\infty(\R^2_+)}\|\nabla\omega(t)\|_{L^2(\R^2_+)}.
\end{eqnarray*}
Therefore, by the Gronwall's inequality, we obtain, for all $t\in[0,T^\ast)$,
\begin{eqnarray}\label{5.16}
\|\nabla\omega(t)\|_{L^2(\R^2_+)}
\leq
\|\nabla\omega_0\|_{L^2(\R^2_+)}e^{\int_{0}^{t}\|\nabla u(\tau)\|_{L^\infty(\R^2_+)}d\tau}.
\end{eqnarray}
By Corollary \ref{single-exponential-velocity-gradient}, we have, for any $t\in[0,T^\ast)$,
\begin{eqnarray}\label{5.17}
\|u(t)\|_{W^{1,\infty}(\R^2_+)}\leq Ce^{Ct},
\end{eqnarray}
where the constants $C=C\Big(p,\|\omega_0\|_{L^\infty(\R^2_+)},\|\omega_0\|_{L^p(\R^2_+)},\|\nabla\omega_0\|_{L^\infty(\R^2_+)}\Big)$.
Combining \eqref{5.16} and \eqref{5.17} yields that, for all $t\in[0,T^\ast)$,
\begin{eqnarray}\label{5.18}
\|\nabla\omega(t)\|_{L^2(\R^2_+)}\leq
\|\nabla\omega_0\|_{L^2(\R^2_+)}e^{Ce^{Ct}}.
\end{eqnarray}
Substituting \eqref{5.17}-\eqref{5.18} into \eqref{5.15} implies
\begin{eqnarray*}
\|\omega(t)\|_{W^{k,p}(\R^2_+)}
\leq \|\omega_0\|_{W^{k,p}(\R^2_+)}e^{Ce^{Ct}}e^{C\|\nabla\omega_0\|_{L^2(\R^2_+)}\int_{0}^{t}e^{Ce^{C\tau}}d\tau},
\end{eqnarray*}
for all $t\in[0,T^\ast)$.
The proof of Theorem \ref{the-2} is complete.
\appendix
\section{EXISTENCE AND UNIQUENESS OF SOLUTIONS TO A LINEAR TRANSPORT EQUATION}
In this appendix, we give the proof of Lemma \ref{Linear-Transport}, which asserts the existence and uniqueness of solutions to a linear transport equation \eqref{4.2}. For convenience, we first state a global approximation by functions smooth up to boundary (see, e.g., \cite{[GT]}).
\begin{lemma}\label{Global-Approximation}
$C^\infty(\bar{\R}^N_+)$ is dense in $W^{m,s}(\R^N_+)$ for $m\geq1$ an integer and $1\leq s<\infty$. That is, for $g\in W^{m,s}(\R^N_+)$, there exists a sequence $\{g^n\}_{n\in\mathbb{N}}$ in $C^\infty(\bar{\R}^N_+)$ such that $g^n$ converges to $g$ in $W^{m,s}(\R^N_+)$.
\end{lemma}
\textbf{Proof}. For $g\in W^{m,s}(\R^N_+)$, we define the approximate functions as the translated mollifications of $g$, given by
\begin{eqnarray*}
g^n(x)=n^N\int_{\mathbb{R}^N_+}\rho\Big(n(x-y)+2e_N\Big)g(y)dy, ~x\in\R^N_+,
\end{eqnarray*}
where $e_N$ denotes the unit coordinate vector in the $x_N$ direction and $\rho$ is a non-negative function in $C_c^\infty(\R^N)$ vanishing outside the unit ball $B_1(0)$ and satisfying $\int_{\mathbb{R}^N}\rho(x)dx=1$.
Then it is readily verified that $g^n$ is in $C^\infty(\bar{\R}^N_+)$ and $g^n$ converges to $g$ in $W^{m,s}(\R^N_+)$. The proof of the lemma is finished.

We are now in a position to prove Lemma \ref{Linear-Transport}.

\textbf{Proof of Lemma \ref{Linear-Transport}.} We will focus on the existence part. The uniqueness part is a direct result of the standard $L^p(\R^2_+)$ energy estimate.
 For any fixed time $t\in[0,T]$, $\omega(t)\in W^{k,p}(\R^2_+)$,
by Lemma \ref{Global-Approximation}, there exist functions $\omega^n(t)\in C^\infty(\bar{\R}^2_+)$ such that
\begin{eqnarray*}
\|\omega^n-\omega\|_{L^\infty (0,T;W^{k,p}(\R^2_+))}\rightarrow0, ~{\rm as}~n\rightarrow\infty,
\end{eqnarray*}
and
\begin{eqnarray*}
\|\omega^n-\omega\|_{C ([0,T];L^p(\R^2_+))}\rightarrow0, ~{\rm as}~n\rightarrow\infty.
\end{eqnarray*}
Therefore $\{\omega^n\}_{n\in\mathbb{N}}$ is uniformly bounded in
$L^\infty(0,T;W^{k,p}(\R^2_+))$.
Define the stream function sequence $\{\Psi^n\}_{n\in\mathbb{N}}$ as follows:
\begin{eqnarray*}
\Psi^n(x,t)=\frac{1}{2\pi}\int_{\mathbb R^2_+}\Big(\log|x-y|-\log|x-\bar{y}|\Big)\omega^n(y,t)dy.
\end{eqnarray*}
Clearly, $\Psi^n$ satisfies the following Poisson equation with homogeneous Dirichlet boundary condition£º
\begin{equation*}
\left\{\ba
&\Delta\Psi^n(\cdot,t)=\omega^n(\cdot,t), ~x\in \R^{2}_+,\\
&\Psi^n(x,t)=0,~x\in \partial\R^{2}_+. \ea\ \right.
\end{equation*}
The velocity field sequence $u^n(\cdot,t):=\nabla^\bot\Psi^n(\cdot,t)$ can be written as
\begin{eqnarray}\label{A.1}
u^n(x,t)=\frac{1}{2\pi}\int_{\mathbb R^2_+}\Big(\frac{(x-y)^\bot}{|x-y|^2}-\frac{(x-\bar{y})^\bot}{|x-\bar{y}|^2}\Big)\omega^n(y,t)dy.
\end{eqnarray}
Clearly, $u^n(x,t)$  is divergence free and tangential to the boundary $\partial\R^{2}_+$,
that is,
\begin{eqnarray}\label{A.2}
{\rm div}~u^n(x,t)=0,~x\in\R^2_+~{\rm and}~ u^n_2(x_1,0,t)=0, ~x_1\in\R.
\end{eqnarray}
By the Sobolev embedding $W^{1,\frac{2p}{2-p}}(\R^{2}_+)\hookrightarrow L^\infty(\R^2_+)$ for $1<p<2$
and Proposition \ref{Sobolev-type estimate}, we have
\begin{eqnarray*}
\|\nabla u^n\|_{L^\infty(0,T;L^\infty(\R^2_+))}
&\leq& C\|\nabla u^n\|_{L^\infty(0,T;W^{1,\frac{2p}{2-p}}(\R^{2}_+))}\\
&\leq& C\|\omega^n\|_{L^\infty(0,T;W^{k,p}(\R^2_+))}.
\end{eqnarray*}
Then the standard Cauchy-Lipchitz theorem ensures that we can define  a smooth flow map $\Phi^n_t(\cdot)$  on $[0,T]\times\bar{\R}^2_+$ through $u^n$.
From the characteristic method, the function $\theta^n:[0,T]\times\bar{\R}^2_+\rightarrow\R$ defined by
\begin{eqnarray*}
\theta^n(x,t)=\omega^n(\Phi^n_{-t}(x),0),
\end{eqnarray*}
is a unique smooth solution to the following linear transport equation
\begin{equation}\label{A.3}
\left\{\ba
&\partial_t\theta^n+u^n\cdot\nabla\theta^n=0, ~(x,t)\in \R^{2}_+\times\R_{+},\\
&\theta^n(\cdot,0)=\omega^n(\cdot,0).\ea\ \right.
\end{equation}
Here $\Phi^n_{-t}$ is the inverse of $\Phi^n_t$.

Next, we  derive some uniform estimates for the approximate sequence $\{\theta^n\}_{n\in\mathbb{N}}$.
It is direct to obtain
\begin{eqnarray}\label{A.4}
\frac{d}{dt}\|\theta^n(t)\|_{L^{p}(\R^2_+)}=0.
\end{eqnarray}
For arbitrary multi-indice $\alpha\in\mathbb{N}^{2}$ with $0<|\alpha|\leq k$, operating $\partial^{\alpha}$ on the first eqaution of \eqref{A.3} gives
\begin{eqnarray}\label{A.5}
\partial_t\partial^{\alpha}\theta^n+u^n\cdot\nabla \partial^{\alpha}\theta^n+\partial^{\alpha}(u^n\cdot\nabla\theta^n)-u^n\cdot\nabla\partial^{\alpha}\theta^n=0.
\end{eqnarray}
Multiplying \eqref{A.5} by $|\partial^{\alpha}\theta^n|^{p-2}\partial^{\alpha}\theta^n$, integrating the resulting equation on $\R^{2}_+$, making integration by parts  and utilizing \eqref{A.2}, H\"{o}lder inequality, Corollary \ref{Calculus-half-space}, \eqref{4.8}, Sobolev embedding $W^{1,\frac{2p}{2-p}}(\R^{2}_+)\hookrightarrow L^\infty(\R^2_+)$ for $1<p<2$ and Proposition \ref{Sobolev-type estimate}, we can finally obtain
\begin{eqnarray*}
&&\frac{1}{p}\frac{d}{dt}\|\partial^{\alpha}\theta^n(t)\|^{p}_{L^{p}(\R^2_+)}\\
&=&-\int_{\mathbb{R}^2_+}\Big(\partial^{\alpha}(u^n\cdot\nabla\theta^n)-u^n\cdot\nabla \partial^{\alpha}\theta^n\Big)|\partial^{\alpha}\theta^n|^{p-2}(\partial^{\alpha}\theta^n)dx\\
&\leq&\|\partial^{\alpha}(u^n\cdot\nabla\theta^n)-u^n\cdot\nabla \partial^{\alpha}\theta^n\|_{L^{p}(\R^2_+)}\|\partial^{\alpha}\theta^n\|^{p-1}_{L^{p}(\R^2_+)}\\
&\leq&C\Big(\|u^n\|_{W^{|\alpha|,\frac{2p}{2-p}}(\R^2_+)}\|\nabla\theta^n\|_{L^2(\R^2_+)}\\
&&\
+\|u^n\|_{W^{1,\infty}(\R^2_+)}\|\nabla\theta^n\|_{W^{|\alpha|-1,p}(\R^2_+)}\Big)\|\partial^{\alpha}\theta^n\|^{p-1}_{L^{p}(\R^2_+)}\\
&\leq&C\|u^n\|_{W^{k,\frac{2p}{2-p}}(\R^2_+)}\|\theta^n\|_{W^{k,p}(\R^2_+)}\|\partial^{\alpha}\theta^n\|^{p-1}_{L^{p}(\R^2_+)}\\
&\leq&C\|\omega^n(t)\|_{W^{k,p}(\R^2_+)}\|\theta^n(t)\|_{W^{k,p}(\R^2_+)}\|\partial^{\alpha}\theta^n\|^{p-1}_{L^{p}(\R^2_+)},
\end{eqnarray*}
where the constant $C$ is independent of $n$.
It follows that
\begin{eqnarray}\label{A.6}
\frac{d}{dt}\|\partial^{\alpha}\theta^n(t)\|_{L^{p}(\R^2_+)}
\leq C\|\omega^n(t)\|_{W^{k,p}(\R^2_+)}\|\theta^n(t)\|_{W^{k,p}(\R^2_+)}
\end{eqnarray}
 for any $0<|\alpha|\leq k$.
Combining \eqref{A.4} with \eqref{A.6} leads to
\begin{eqnarray*}
\frac{d}{dt}\|\theta^n(t)\|_{W^{k,p}(\R^2_+)}
\leq C\|\omega^n(t)\|_{W^{k,p}(\R^2_+)}\|\theta^n(t)\|_{W^{k,p}(\R^2_+)}.
\end{eqnarray*}
Using the Gronwall's inequality, we  get
\begin{eqnarray}\label{A.7}
\|\theta^n(t)\|_{W^{k,p}(\R^2_+)}
&\leq& \|\theta^n(0)\|_{W^{k,p}(\R^2_+)}e^{C{\int_{0}^{t}\|\omega^n(\tau)\|_{W^{k,p}(\R^2_+)}d\tau}}\\
&\leq& \|\omega^n(0)\|_{W^{k,p}(\R^2_+)}e^{CT\|\omega^n\|_{L^\infty(0,T;W^{k,p}(\R^2_+))}}.\nonumber
\end{eqnarray}
Since $\{\omega^n\}_{n\in\mathbb{N}}$ is uniformly bounded in
$L^\infty(0,T;W^{k,p}(\R^2_+))$, it follows that $\{\theta^n\}_{n\in\mathbb{N}}$ is uniformly bounded in $L^\infty(0,T;W^{k,p}(\R^2_+))$.
Furthermore, utilizing \eqref{A.3}, Corollary \ref{Calculus-half-space}, Proposition \ref{Sobolev-type estimate}, \eqref{4.8} and the Sobolev embedding $W^{1,\frac{2p}{2-p}}(\R^{2}_+)\hookrightarrow L^\infty(\R^2_+)$ for $1<p<2$, we can arrive at
\begin{eqnarray*}
&&\|\partial_t\theta^n(t)\|_{W^{k-1,p}(\R^2_+)}\\
&=& \|u^n\cdot\nabla\theta^n\|_{W^{k-1,p}(\R^2_+)}\\
&\leq&C\Big(\|u^n\|_{W^{k-1,\frac{2p}{2-p}}(\R^2_+)}\|\nabla\theta^n\|_{L^2(\R^2_+)}
+\|u^n\|_{W^{1,\infty}(\R^2_+)}\|\nabla\theta^n\|_{W^{k-1,p}(\R^2_+)}\Big) \\
&\leq&C\|\omega^n(t)\|_{W^{k,p}(\R^2_+)}\|\theta^n(t)\|_{W^{k,p}(\R^2_+)},
\end{eqnarray*}
 which shows that $\{\partial_t\theta^n(t)\}_{n\in\mathbb{N}}$ is uniformly bounded in $L^\infty(0,T;W^{k-1,p}(\R^2_+))$.
 Then $\{\theta^n\}_{n\in\mathbb{N}}$ is also uniformly bounded in ${\rm Lip}([0,T];W^{k-1,p}(\R^2_+))$.
 Next we show that $\{\theta^n\}_{n\in\mathbb{N}}$ is a Cauchy sequence in $C([0,T];L^{p}(\R^2_+))$.
 It follows from \eqref{A.3} that
\begin{equation}\label{A.8}
\left\{\ba
&\partial_t\theta^m+u^m\cdot\nabla\theta^m=0,\\
&\partial_t\theta^n+u^n\cdot\nabla\theta^n=0. \ea\ \right.
\end{equation}
Subtracting the equations in \eqref{A.8} gives
\begin{eqnarray*}
\partial_t(\theta^m-\theta^n)+u^m\cdot\nabla(\theta^m-\theta^n)+(u^m-u^n)\cdot\nabla\theta^n=0.
\end{eqnarray*}
Then, similar to \eqref{4.7}, we can obtain
\begin{eqnarray*}
\frac{d}{dt}\|\theta^m(t)-\theta^n(t)\|_{L^{p}(\R^2_+)}
\leq C\|\theta^n(t)\|_{W^{k,p}(\R^2_+)}\|\omega^m(t)-\omega^n(t)\|_{L^{p}(\R^2_+)},
\end{eqnarray*}
which shows that $\{\theta^n\}_{n\in\mathbb{N}}$ is a Cauchy sequence in $C([0,T];L^{p}(\R^2_+))$.
 By the Gagliardo-Nirenberg inequality
$\|f\|_{W^{k-1,p}(\R^2_+)}\leq C\|f\|^{\frac{1}{k}}_{L^p(\R^2_+)}\|f\|^{1-\frac{1}{k}}_{W^{k,p}(\R^2_+)} $
and the uniform boundness of $\{\theta^n\}_{n\in\mathbb{N}}$ in $L^\infty(0,T;W^{k,p}(\R^2_+))$ , we can find that the sequence
 $\{\theta^n\}_{n\in\mathbb{N}}$ is also Cauchy in $C([0,T];W^{k-1,p}(\R^2_+))$.
Thus, there exists a $\theta$ in $L^\infty(0,T;W^{k,p}(\R^2_+))\cap{\rm Lip}([0,T\rfloor;W^{k-1,p}(\R^2_+))$ such that
\begin{eqnarray}\label{A.9}
\theta_n\rightarrow\theta~{\rm in}~C([0,T];W^{k-1,p}(\R^2_+)),~{\rm as}~n\rightarrow\infty,
\end{eqnarray}
\begin{eqnarray}\label{A.10}
\theta_n\overset{\ast}{\rightharpoonup}\theta~{\rm in}~L^\infty(0,T;W^{k,p}(\R^2_+)),~{\rm as}~n\rightarrow\infty,
\end{eqnarray}
and
\begin{eqnarray}\label{A.11}
\partial_t\theta_n\overset{\ast}{\rightharpoonup}\partial_t\theta~{\rm in}~L^\infty(0,T;W^{k-1,p}(\R^2_+)),~{\rm as}~n\rightarrow\infty.
\end{eqnarray}
Moreover, by Proposition \ref{Sobolev-type estimate}, we have $\|u^n-u\|_{W^{k,\frac{2p}{2-p}}(\R^2_+)}\leq C_{k,p}\|\omega^n-\omega\|_{W^{k,p}(\R^2_+)}$. Thus $u^n$ converges to $u$ in $L^\infty(0,T;W^{k,\frac{2p}{2-p}}(\R^2_+))$. Then it follows from \eqref{A.9}-\eqref{A.11} that $\theta$ is a solution to \eqref{4.2}.
Also, from \eqref{A.7}, we can obtain, for any $t\in[0,T]$,
\begin{eqnarray*}
\|\theta(t)\|_{W^{k,p}(\R^2_+)}
\leq\|\omega_0\|_{W^{k,p}(\R^2_+)}
e^{C{\int_{0}^{t}\|\omega(\tau)\|_{W^{k,p}(\R^2_+)}d\tau}}.
\end{eqnarray*}
The proof of the lemma is complete.

In the end of this appendix, we present the expression of $\nabla u$ in the half plane for completeness. Similar computations for the whole plane case can be found in \cite{[MB]}.
\begin{lemma}\label{Formula for gradient}
Let the velocity field $u(x,t)$ be  defined by
\begin{eqnarray*}
u(x,t)=\int_{\mathbb R^2_+}K(x,y)\omega(y,t)dy,
\end{eqnarray*}
where $K(x,y)=\frac{1}{2\pi}\Big(\frac{(x-y)^\bot}{|x-y|^2}-\frac{(x-\bar{y})^\bot}{|x-\bar{y}|^2}\Big)$, and $\omega$ is the vorticity of $u$.
Then
\begin{eqnarray*}
\nabla u(x,t)=P.V.\int_{\R_+^2}\nabla_x K(x,y)\omega(y,t)dy+\frac{\omega(x,t)}{2}
\begin{pmatrix}
0&1\\-1&0
\end{pmatrix}.
\end{eqnarray*}
Furthermore,
\begin{eqnarray*}
\nabla u(x,t)=\frac{1}{2\pi}P.V.\int_{\R^2_+}M(x,y)\omega(y,t)dy+\frac{\omega(x,t)}{2}
\begin{pmatrix}
0&1\\-1&0
\end{pmatrix},
\end{eqnarray*}
where the kernel matrix
\begin{eqnarray*}
M(x,y)=
\begin{pmatrix}
\frac{-2(x_1-y_1)(x_2-y_2)}{|x-y|^4}+\frac{2(x_1-y_1)(x_2+y_2)}{|x-\bar{y}|^4}
&\frac{(x_1-y_1)^2-(x_2-y_2)^2}{|x-y|^4}-\frac{(x_1-y_1)^2-(x_2+y_2)^2}{|x-\bar{y}|^4}\\
\frac{(x_1-y_1)^2-(x_2-y_2)^2}{|x-y|^4}-\frac{(x_1-y_1)^2-(x_2+y_2)^2}{|x-\bar{y}|^4}
&\frac{2(x_1-y_1)(x_2-y_2)}{|x-y|^4}-\frac{2(x_1-y_1)(x_2+y_2)}{|x-\bar{y}|^4}
\end{pmatrix}.
\end{eqnarray*}
Consequently, the velocity field $u$ satisfies the divergence-free condition
\begin{eqnarray*}
{\rm div}~u(x,t)=0,~x\in\R^2_+.
\end{eqnarray*}
\end{lemma}

\textbf{Proof}.
For any $1\leq i,j\leq2$ and $\varphi\in C_c^\infty(\R^2_+)$, by the definition of the distributional derivative and Fubini's theorem, we have
\begin{eqnarray}\label{A.12}
\langle\partial_ju^i,\varphi\rangle
&=&-\langle u^i,\partial_j\varphi\rangle\nonumber\\
&=&-\int_{\mathbb R^2_+}u^i(x)\partial_{x_j}\varphi(x)dx\nonumber\\
&=&-\int_{\mathbb R^2_+}\Big(\int_{\mathbb R^2_+}K^i(x,y)\omega(y,t)dy\Big)\partial_{x_j}\varphi(x)dx\nonumber\\
&=&-\int_{\mathbb R^2_+}\omega(y,t)\Big(\int_{\mathbb R^2_+}K^i(x,y)\partial_{x_j}\varphi(x)dx\Big)dy\nonumber\\
&=&-\int_{\mathbb R^2_+}\omega(y,t)I_{ij}(y)dy,
\end{eqnarray}
where the integrals $I_{ij}(y):=\int_{\mathbb R^2_+}K^i(x,y)\partial_{x_j}\varphi(x)dx$.

By the dominated convergence theorem and integration by parts, $I_{ij}(y)$ is expressed as
\begin{eqnarray}\label{A.13}
I_{ij}(y)
&=&\int_{\mathbb R^2_+}K^i(x,y)\partial_{x_j}\varphi(x)dx\nonumber\\
&=&\displaystyle\lim_{\varepsilon\downarrow0^+}\int_{\mathbb{R}^2_+\cap(B_{\varepsilon}(y))^c}K^i(x,y)\partial_{x_j}\varphi(x)dx\nonumber\\
&=&
\displaystyle\lim_{\varepsilon\downarrow0^+}\Big(\int_{\partial\mathbb{R}^2_+}K^i(x,y)\varphi(x)n_j(x)dS_x
-\int_{|x-y|=\varepsilon}K^i(x,y)\varphi(x)\frac{x_j-y_j}{|x-y|}dS_x\Big)\nonumber\\
&&\ -\displaystyle\lim_{\varepsilon\downarrow0^+}\int_{\mathbb{R}^2_+\cap(B_{\varepsilon}(y))^c}\partial_{x_j}K^i(x,y)\varphi(x)dx\nonumber\\
&=&-\displaystyle\lim_{\varepsilon\downarrow0^+}\int_{\mathbb{R}^2_+\cap(B_{\varepsilon}(y))^c}\partial_{x_j}K^i(x,y)\varphi(x)dx\nonumber\\
&&\
-\displaystyle\lim_{\varepsilon\downarrow0^+}\varepsilon\int_{|z|=1}K^i(y+\varepsilon z,y)\varphi(y+\varepsilon z)z_jdS_z.
\end{eqnarray}
Therefore, using the Fubini's theorem again and \eqref{A.12}-\eqref{A.13}, we can obtain
\begin{eqnarray}\label{A.14}
\langle\partial_ju^i,\varphi\rangle
&=&\int_{\mathbb R^2_+}\omega(y,t)\Big(\displaystyle\lim_{\varepsilon\downarrow0^+}\int_{\mathbb{R}^2_+\cap(B_{\varepsilon}(y))^c}\partial_{x_j}K^i(x,y)\varphi(x)dx\Big) dy\nonumber\\
&&\
+\int_{\mathbb R^2_+}\omega(y,t)\Big(\displaystyle\lim_{\varepsilon\downarrow0^+}\varepsilon\int_{|z|=1}K^i(y+\varepsilon z,y)\varphi(y+\varepsilon z)z_jdS_z\Big)dy\nonumber\\
&=&\int_{\mathbb R^2_+}\Big(\displaystyle\lim_{\varepsilon\downarrow0^+}\int_{\mathbb{R}^2_+\cap(B_{\varepsilon}(x))^c}\partial_{x_j}K^i(x,y)\omega(y,t)dy\Big)\varphi(x)dx \nonumber\\
&&\
+\int_{\mathbb R^2_+}\omega(y,t)\Big(\displaystyle\lim_{\varepsilon\downarrow0^+}\varepsilon\int_{|z|=1}K^i(y+\varepsilon z,y)\varphi(y+\varepsilon z)z_jdS_z\Big)dy\nonumber\\
&=&\int_{\mathbb R^2_+}\Big(P.V.\int_{\mathbb{R}^2_+}\partial_{x_j}K^i(x,y)\omega(y,t)dy\Big)\varphi(x)dx \nonumber\\
&&\
+\int_{\mathbb R^2_+}\omega(y,t)\Big(\displaystyle\lim_{\varepsilon\downarrow0^+}\varepsilon\int_{|z|=1}K^i(y+\varepsilon z,y)\varphi(y+\varepsilon z)z_jdS_z\Big) dy\nonumber\\
&=&\langle P.V.\int_{\mathbb{R}^2_+}\partial_{x_j}K^i(x,y)\omega(y,t)dy,\varphi\rangle\nonumber\\
&&\
+\int_{\mathbb R^2_+}\omega(y,t)\Big(\displaystyle\lim_{\varepsilon\downarrow0^+}\varepsilon\int_{|z|=1}K^i(y+\varepsilon z,y)\varphi(y+\varepsilon z)z_jdS_z\Big) dy.
\end{eqnarray}
The dominated convergence theorem and the identity $\int_{|z|=1}z_jz_2dS_z=\delta_{2j}\pi$ deduce that
\begin{eqnarray}\label{A.15}
&&\displaystyle\lim_{\varepsilon\downarrow0^+}\varepsilon\int_{|z|=1}K^1(y+\varepsilon z,y)\varphi(y+\varepsilon z)z_jdS_z\nonumber\\
&=& \frac{1}{2\pi}\displaystyle\lim_{\varepsilon\downarrow0^+}\varepsilon\int_{|z|=1}\Big(\frac{y_2+\varepsilon z_2-y_2}{|y+\varepsilon z-y|^2}-\frac{y_2+\varepsilon z_2+y_2}{|y+\varepsilon z-\bar{y}|^2}\Big)\varphi(y+\varepsilon z)z_jdS_z\nonumber\\
&=&\frac{1}{2\pi}\displaystyle\lim_{\varepsilon\downarrow0^+}\int_{|z|=1}\varphi(y+\varepsilon z)z_jz_2dS_z
-\frac{1}{2\pi}\displaystyle\lim_{\varepsilon\downarrow0^+}\int_{|z|=1}\frac{\varepsilon(2y_2+\varepsilon z_2)}{|2y_2e_2+\varepsilon z|^2}\varphi(y+\varepsilon z)z_jdS_z\nonumber\\
&=&\frac{1}{2\pi}\varphi(y)\int_{|z|=1}z_jz_2dS_z\nonumber\\
&=&\frac{1}{2}\varphi(y)\delta_{2j},
\end{eqnarray}
where $\delta_{ij}$ is the Kronecker symbol defined by
\begin{eqnarray*}
\delta_{ij}=
\begin{cases}
1,
& \mbox{if $i=j,$ } \\
0,
& \mbox{if $i\neq j.$ } \\
\end{cases}
\end{eqnarray*}
Similarly, we have
\begin{eqnarray}\label{A.16}
&&\displaystyle\lim_{\varepsilon\downarrow0^+}\varepsilon\int_{|z|=1}K^2(y+\varepsilon z,y)\varphi(y+\varepsilon z)z_jdS_z\nonumber\\
&=& \frac{1}{2\pi}\displaystyle\lim_{\varepsilon\downarrow0^+}\varepsilon\int_{|z|=1}\Big(\frac{y_1-(y_1+\varepsilon z_1)}{|y+\varepsilon z-y|^2}-\frac{y_1-(y_1+\varepsilon z_1)}{|y+\varepsilon z-\bar{y}|^2}\Big)\varphi(y+\varepsilon z)z_jdS_z\nonumber\\
&=&-\frac{1}{2\pi}\displaystyle\lim_{\varepsilon\downarrow0^+}\int_{|z|=1}\varphi(y+\varepsilon z)z_1z_jdS_z
+\frac{1}{2\pi}\displaystyle\lim_{\varepsilon\downarrow0^+}\int_{|z|=1}\frac{\varepsilon^2z_1}{|2y_2e_2+\varepsilon z|^2}\varphi(y+\varepsilon z)z_jdS_z\nonumber\\
&=&-\frac{1}{2\pi}\varphi(y)\int_{|z|=1}z_1z_jdS_z\nonumber\\
&=&-\frac{1}{2}\varphi(y)\delta_{1j},
\end{eqnarray}
where we have used $\int_{|z|=1}z_1z_jdS_z=\delta_{1j}\pi$.

Combining  \eqref{A.14}-\eqref{A.16}, it then holds that, for every $1\leq j\leq2$,
\begin{eqnarray}\label{A.17}
\langle\partial_ju^1,\varphi\rangle
&=&\langle P.V.\int_{\mathbb{R}^2_+}\partial_{x_j}K^1(x,y)\omega(y,t)dy,\varphi\rangle
+\int_{\mathbb R^2_+}\omega(y,t)\frac{1}{2}\varphi(y)\delta_{2j}dy\nonumber\\
&=&\langle P.V.\int_{\mathbb{R}^2_+}\partial_{x_j}K^i(x,y)\omega(y,t)dy+\frac{1}{2}\delta_{2j}\omega(\cdot,t),\varphi\rangle,
\end{eqnarray}
and
\begin{eqnarray}\label{A.18}
\langle\partial_ju^2,\varphi\rangle
&=&\langle P.V.\int_{\mathbb{R}^2_+}\partial_{x_j}K^2(x,y)\omega(y,t)dy,\varphi\rangle
-\int_{\mathbb R^2_+}\omega(y,t)\frac{1}{2}\varphi(y)\delta_{1j}dy\nonumber\\
&=&\langle P.V.\int_{\mathbb{R}^2_+}\partial_{x_j}K^i(x,y)\omega(y,t)dy-\frac{1}{2}\delta_{1j}\omega(\cdot,t),\varphi\rangle.
\end{eqnarray}
By the definition of the distribution derivatives and \eqref{A.17}-\eqref{A.18}, we conclude the proof of the first part of Lemma \ref{Formula for gradient}.
The rest of the proof can be deduced from some standard computations.

{\bf Acknowledgements.}
Q. Jiu was partially supported by the National Natural Science Foundation of
China (NNSFC) (No. 11931010, No. 12061003)
and key research project of the Academy for Multidisciplinary Studies of CNU.
Y. Li was partially supported
by Beijing Municipal Natural Science Foundation (1212002).

\end{document}